%% file: paper.tex
\documentclass[11pt]{amsart}
\usepackage{amsmath}
\usepackage{geometry}
\usepackage{adjustbox}
\usepackage{multirow}
\usepackage{cases}
\usepackage[hidelinks]{hyperref}
\usepackage{xcolor}
\usepackage{hyperref}
\usepackage{cleveref}
\usepackage{algorithm}
\usepackage{forest}
\usepackage{tikz}
\usetikzlibrary{positioning,quotes,spy,shapes,arrows.meta,arrows,chains,calc,arrows,decorations.markings}
\usepackage{algpseudocode}
\usepackage{mathtools}
\usetikzlibrary{trees}
\title{Neural-network-based regularization methods for inverse problems in imaging}
\author{Andreas Habring \and Martin Holler}

\usepackage{booktabs}
\include{notation}

\pgfdeclarelayer{bg}    %
\pgfsetlayers{bg,main}
\definecolor{UniGrazWhite}{RGB}{255,255,255}
\definecolor{UniGrazYellow}{RGB}{255,213,0}
\definecolor{UniGrazGray}{RGB}{204,204,204}
\definecolor{UniGrazBlack}{RGB}{0,0,0}
\definecolor{UniGrazNAWI}{RGB}{6,120,180}
\definecolor{UniGrazMath}{RGB}{144,174,181}
\definecolor{UniGrazBlue}{RGB}{31,31,139}
\definecolor{green}{rgb}{0.0, 0.5, 0.0}
\definecolor{orange}{rgb}{1.0, 0.44, 0.37}

\usetikzlibrary{spy,decorations.pathreplacing,angles,quotes,calc,arrows,positioning}

\definecolor{spycolor}{RGB}{150,150,200}
\tikzset{
    rectspy/.default={lens={scale=3}, size=3cm},
    rectspy on/.style={#1,},
    rectspy/.style={
        draw=spycolor,
        connect spies,
        spy scope={
        every spy on node/.style={
            draw=spycolor,
            very thick,
            rectangle, 
            rectspy on,
        },
        every spy in node/.style={
            draw=spycolor,
            very thick,
            rectangle,
        },
        #1
        },
        spy connection path={\draw[spycolor, very thick] (tikzspyonnode) -- (tikzspyinnode);}
    }
}

\tikzset{
    sepbar/.style={
        very thick,
        black!30!white,
    }
}

\begin{document}

\begin{abstract}
This review provides an introduction to - and overview of - the current state of the art in neural-network based regularization methods for inverse problems in imaging. It aims to introduce readers with a solid knowledge in applied mathematics and a basic understanding of neural networks to different concepts of applying neural networks for regularizing inverse problems in imaging. Distinguishing features of this review are, among others, an easily accessible introduction to learned generators and learned priors, in particular diffusion models, for inverse problems, and a section focusing explicitly on existing results in function space analysis of neural-network-based approaches in this context.
\end{abstract}

\maketitle

\setcounter{tocdepth}{1}
\tableofcontents

\begin{figure}
\begin{tikzpicture}[inner sep=0,rectspy={lens={scale=3}, width=4.9cm, height=1.7cm}]
\node at (0,0) {\includegraphics[width=4.9cm]{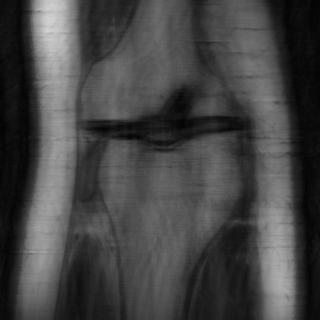}};
\spy on (1.3,0.65) in node at (0,-3.4);
\node at (0,-4.7) {Zero-fill reconstruction};

\node at (5.2,0) {\includegraphics[width=4.9cm]{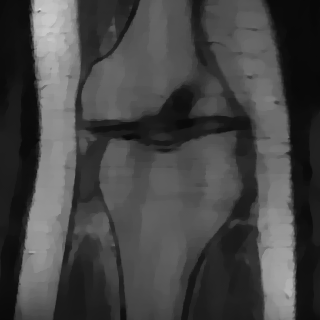}};
\spy on (6.5,0.65) in node at (5.2,-3.4);
\node at (5.2,-4.7){Hand crafted TV prior};

\node at (10.4,0) {\includegraphics[width=4.9cm]{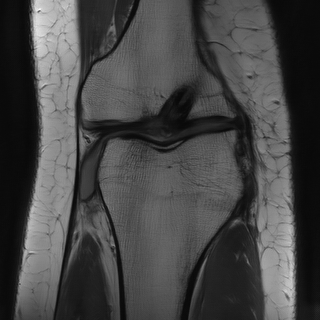}};
\spy on (11.7,0.65) in node at (10.4,-3.4);
\node at (10.4,-4.7) {CNN-based reconstruction};

\end{tikzpicture}
\caption{\label{fig:mri_example} Comparison of different methods for the reconstruction of magnetic resonance images from subsampled data. From left to right: Reconstruction by setting unknown data components to zero, classical variational approach with TV prior \cite{rudin1992nonlinear}, result of trained, neural-network based reconstruction approach from \cite{hammernik2018learning}, see also \cite{zach2023stable}$^1$. Data taken from \cite{knoll2020fastmri}.}
\end{figure}
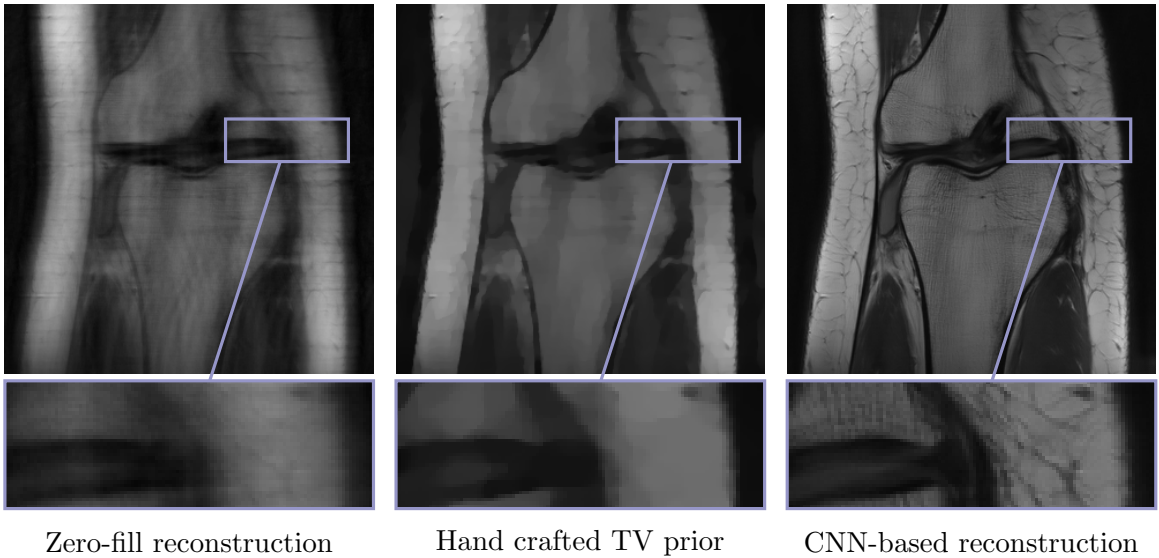

\section{Introduction}\label{sec:introduction}
The field of inverse problems in imaging is probably among the mathematical fields that were most affected by the rise of neural-network-based deep learning techniques of the past decade. 
Accounting for this, the present article strives to provide an overview of the current state of the art regarding techniques to solve inverse problems in imaging under the use of neural networks (NNs), in particular convolutional neural networks (CNNs). Its main focus is to explain the most successful approaches conceptually and, thus, enable a clear unterstanding of the discussed methods rather than provide an exhaustive literature review. In order to facilitate a concise presentation, details regarding notation and precise mathematical definitions are kept to a minimum, and readers are assumed to have a good basis in applied mathematics and to be already familiar with the notion of neural networks in general. Furthermore, we emphasize already at this point that, given the size and rapid development of the field, such a review can of course never be complete and we refer readers also to other reviews listed in \Cref{sec:meta_review} below for a more comprehensive overview.

In inverse problems in imaging, the goal is to infer an unknown ground truth signal $u\in\Xc$ (usually an image) given a measurement or observation $y\in\Yc$. Ground truth and observation are related to each other via the so-called forward operator $\fwd:\Xc\rightarrow\Yc$ according to
\begin{equation}\label{eq:inverse_problem}
y = \fwd  u + \zeta,
\end{equation}\footnotetext[1]{The authors thank Martin Zach from Graz University of Technology for providing these results.}
where $\zeta$ denotes possible measurement noise (we consider additive noise here for simplicity, but the noise might, e,g., also be multiplicative). In many cases of interest such as magnetic resonance imaging (MRI), computed tomography (CT), electrical impedance tomography (EIT) or super-resolution the forward operator is, in fact, not invertible or its inverse is (severely) ill-conditioned, rendering \eqref{eq:inverse_problem} an ill-posed problem. A popular remedy to overcome ill-posedness are variational regularization methods \cite[Section 3]{scherzer2008variational} where, as a surrogate for \eqref{eq:inverse_problem}, one aims to solve
\begin{equation}\label{eq:var_reg}
\min\limits_{u}\Dc_y(\fwd u) + \lambda\Rc(u)
\end{equation}
where $\Dc_y(\fwd u)$ denotes a data fidelity term ensuring $\fwd u\approx y$, such as, e.g., $\Dc_y(\fwd u) = \frac{1}{2}\|\fwd u -y\|_2^2$. The functional $\Rc$ denotes a regularizer which ensures that \eqref{eq:var_reg} is well-posed, i.e., it admits a solution which depends continuously on the data \cite[Section 3.1]{scherzer2008variational}. An example highly popular in particular in imaging problems is the total variation functional (TV) \cite{rudin1992nonlinear,bvfunctions}, which reads as $\Rc(u) = \TV(u) = \|\nabla u\|_1$ in the discrete setting. 

Especially in the context of data-driven approaches, a Bayesian perspective of \eqref{eq:inverse_problem} is often adopted. In this case $u,y,\zeta$ are modeled as random variables, with $u,y$ having a joiont distribution $p(u,y)$, and one considers the log of the \emph{posterior} distribution of $u$ after observing $y$, which can be computed using Bayes' theorem as
\begin{equation}\label{eq:posterior}
-\log \p(u|y) = -\log\p(y|u) - \log \p(u) + \log \p(y).
\end{equation}
We call $\p(y|u)$ the data likelihood and $\p(y)$ the prior distribution. As is common in machine learning literature, as long as there is no risk of ambiguity, we simply use the letter $\p$ to denote any probability density and specify the exact density only through the input argument, e.g., $\p(u)$ denotes the density of the ground truth $u$. Note that \eqref{eq:posterior} allows for point estimation of $u$ via computing the maximum a-posteriori (MAP) estimate, that is, the value of $u$ maximizing $\p(u|y)$. On the other hand, Markov chain Monte Carlo (MCMC) methods such as the unadjusted Langevin algorithm and its variations \cite{durmus19langevin,durmus2019analysis,durmus2022proximal,laumont2022bayesian,habring2023subgradient} allow for sampling from $\p(u|y)$. Thus, instead of the MAP, we can as well compute the expected value of $\p(u|y)$ called the MMSE or, for that matter, any other statistic of the posterior such as the variance \cite{narnhofer2022posterior,zach2022computed,zach2023stable}.

If we identify the distributions occurring in \eqref{eq:posterior} as the Gibbs measures 
\[ \p(y|u) \propto  \exp(-\Dc_y(\fwd u))\quad \text{and}\quad \p(u) \propto \exp(-\Rc(u)), \] we find an equivalence of variational regularization and Bayesian modeling, which is why a strong distinction between the two approaches is not necessary in many cases and we will frequently discuss, e.g., models for $\p(u)$ while keeping in mind that this is equivalent to models of $\Rc(u)$.

While the modeling of $\p(y|u)$ is mostly a question of obtaining appropriate forward and noise models (which is an important research direction on its own), the prior distribution, respectively regularizer, has been subject to a lot of research in particular in inverse imaging problems, with popular examples being Tikhonov regularization \cite[Section 5]{engl1996regularization}, total variation regularization \cite{rudin1992nonlinear} and its generalization to higher-order approaches \cite{bredies2010tgv,holler20ip_review_mh}.
As a consequence of the high complexity of image data, hand crafted methods for solving inverse problems seemed to reach a limit in terms of empirical performance. However, in the past decades, data driven approaches were able to significantly push forward the boundaries of reconstruction quality with appropriate priors. An illustrative example in this context is provided in \Cref{fig:mri_example}, were the performance of the hand crafted TV prior \cite{rudin1992nonlinear} is compared to the performance of the neural-network based approach of \cite{hammernik2018learning} for subsampled MRI.

In this article we focus our attention specifically on data driven methods utilizing NNs. While there exist other regimes, the presented methods focus on the case that the forward operator $\fwd$ is known and will be classified in three different categories: i) unsupervised methods, where during training we require only a sample $(u_i)_i$ of images from the prior distribution $\p(u)$ (\Cref{sec:unsupervised}), ii) supervised methods, where during training we use a sample $(u_i,y_i)_i$ of image-observation pairs from the joint distribution $\p(u,y)$ (\Cref{sec:supervised}), and iii) untrained methods, where NNs are used without any prior training (\Cref{sec:untrained}). This terminology slightly deviates from the conventional one in machine learning literature in the sense that in the context of inverse problems the goal is to predict $u$ from $y$ meaning that $u$, in fact, constitutes the \emph{label} and $y$ the \emph{feature}. Thus, in the classical sense unsupervised would refer to methods training solely on a sample $(y_i)_i$. These approaches are frequently referred to as self-supervised \cite{senouf2019self,hendriksen2020noise,lehtinen2018noise,ongie2020deep} in the context of inverse problems and will not be discussed in detail in the present article. For an overview of the different types of settings, methods and corresponding training strategies considered in this article we refer to \Cref{table:overview}.

After providing an overview of NN based machine learning approaches within these three categories, we further discuss the state of the art in this area in terms of function space modeling (\Cref{sec:function_space}) and uncertainty quantification (\Cref{sec:uncertainty_quantification}).

\subsection{Relation to Existing Review Articles} \label{sec:meta_review} Given the importance of the field, naturally there exist several different, also rather recent review articles on topics related to the scope of this work. Nevertheless, we believe the present article, with a clearly defined scope and a focus on presenting methods, concepts, and mathematical aspects rather then listing existing works, provides an additional contribution. In particular, we believe that some of the unique features of the present article are i) the focus on neural-network based methods for inverse problems in imaging, which allows a comprehensive yet rather detailed, method-driven introduction to the field, ii) a fairly detailed, self-contained summary of the state of the art in using learned generators and learned priors for inverse problems in imaging (see \Cref{sec:learned_generators,sec:learned_prior}), including in particular score based models and diffusion models, and iii) the inclusion of a section focused on existing results in function space analysis (see \Cref{sec:function_space}), which is an important topic in classical inverse problems literature. It is also worth mentioning that the field of deep learning is still evolving rapidly, such that review articles become outdated already after five years or less. This allows the present article to provide a contribution also compared older articles, even if they have a very similar scope.

Nevertheless, many excellent, recent review article exist that cover important aspects that are also in the focus of this work. The following paragraph lists some of them, and we recommend the interested reader to consult them in particular either for a more in-depth treatment of some topics mentioned here only briefly, or to get a larger picture including aspects of neural-network based methods for inverse problems in imaging that are not covered here.

A review article on data-driven methods that is already classical in the mathematical inverse problems community is the comprehensive work \cite{arridge2019solving}. Here, the scope is more general than the one of our work, however, due to the appearance of \cite{arridge2019solving} in 2019, important parts in particular of \Cref{sec:learned_generators,sec:learned_prior} are not covered in the same depth as in the present article.

The scope of \cite{ongie2020deep} is probably closest to the scope of our work. In \cite{ongie2020deep} a taxonomy to categorize different problems and reconstruction methods in the context of deep learning for inverse problems, depending on the extent to which the forward model is known or unknown is presented. Differences to our work are, for instance, that the recency of our work allowed us provide more details on concepts and methods such as score based models, diffusion models and function space theory. 

The very recent review \cite{mukherjee2023learned} focuses in particular on methods with convergence guarantees. It has intersections in particular with \Cref{sec:function_space} of our work, but does not cover topics such as learned generators and learned priors  (see \Cref{sec:learned_generators,sec:learned_prior}) in the same depth as we do.

The review paper \cite{scarlett2022theoretical} is focused on a theoretical perspective of deep learning, in particular in a compressed sensing context, and puts most emphasis on generator-based models and untrained approaches. The work \cite{berner2021modern} in turn provides a review focused on mathematical aspects of deep learning in general, such as approximation results, generalization, and optimization aspects.

A related review is further \cite{shlezinger2023model}, which focuses on model-based deep learning, in particular on the three topics of learning model-based optimization algorithms, unfolding optimizers into trainable architectures (cf. \Cref{sec:unrolling}), and the augmentation of model-based algorithms with trainable neural networks.

Further reviews very much related to the scope of our work are \cite{mccann2017convolutional,lucas2018using}. Due to their appearance in 2017 and 2018, however, they naturally do not fully cover the state of the art in this rapidly evolving field anymore.

Another relevant review is \cite{bai2020deep} which, compared to our work, is probably less method-driven, and again does not cover important parts of the works we discuss in \Cref{sec:learned_generators,sec:learned_prior}.

For a recent review paper on uncertainty quantification in deep learning in general, we refer to \cite{abdar2021review} and for a more extensive review on generator-based regularization for inverse problems see \cite{duff2021regularising}

Finally, the works \cite{li2021review,zhang2020review} provide reviews of deep learning methods in the context of medical imaging.

\section{Notation}

Part of our notation was already introduced in the introduction, and due to the nature of this work, we strive to use notation in a way allowing the educated reader to easily follow our elaborations without the necessity of precise definitions and details.

As already mentioned, the forward model in this article will always be denoted by $\fwd :\Xc \rightarrow \Yc$, where $\Xc$ is the space of images and $\Yc$ is the measurement or observations space, which is not necessarily an image space (cf. MRI). Most of the methods in this article regard $\Xc$ and $\Yc$ to be finite dimensional vectors spaces; the function space settings will only play a role in \Cref{sec:function_space}. Also, for the sake of simplicity and since or focus is on the regularization rather than the forward model, $\fwd$ will always be assumed to be linear unless explicitly stated otherwise. 

Neural networks will usually be denoted as $\Nc_\theta$, where $\theta$ summarizes the network parameters usually consisting of weights (or convolution kernels) $w$ and biases $b$, i.e., $\theta = (w,b)$. Depending on the method at hand, networks will for example map from images space to $\R$, i.e., $\Nc_\theta:\Xc \rightarrow \R$, from latent space to image space, i.e., $\Nc_\theta: \Zc \rightarrow \Xc$, where $\Zc$ and $z$ usually denote the latent space and latent variables, respectively, or from image space to image space $\Nc_\theta: \Xc \rightarrow \Xc $. We assume the reader to be familiar with the general architecture of neural networks, and refer for instance to \cite{berner2021modern} for a precise, mathematical definition.

\begin{table}[h]
\centering
        \begin{tabular}{llll}
        \toprule 
        Setting & Method & Parameter fitting strategy\\
        \toprule
            \multirow{8}{*}[-18pt]{Unsupervised} 	& \multirow{4}{*}[-6.5pt]{Learned generator} 	& Maximum likelihood \\ \cmidrule{3-3} 
            																		& & Autoencoder \\ \cmidrule{3-3} 
            																		& & VAE \\\cmidrule{3-3} 
            																		& & GAN \\\cmidrule{2-3} 
                                 		& \multirow{2}{*}[-2pt]{Learned prior} & Maximum likelihood \\\cmidrule{3-3}
                                 										& & Score-matching and diffusion models\\\cmidrule{2-3} 
                                 		& Plug and Play & Pre-trained denoiser \\\cmidrule{2-3}
                                 		& RED & Pre-trained denoiser \\ \cmidrule{1-3} 
            \multirow{6}{*}[-12pt]{Supervised} & Fully learned & Minimal reconstruction error\\\cmidrule{2-3}
                                 		& \multirow{2}{*}[-2pt]{Post-processing} & Minimal reconstruction error \\\cmidrule{3-3}
                                 										& & GAN \\\cmidrule{2-3}
                                 		& Unrolling & Minimal reconstruction error \\\cmidrule{2-3}
                                 		& \multirow{2}{*}[-2pt]{Learned prior} & Minimal reconstruction error \\\cmidrule{3-3}
                                 										& & Adversarial regularizer \\\cmidrule{1-3}
			Untrained &  Generator-network & Approximate observed data                            										\\
\bottomrule
        \end{tabular}
\vspace*{0.2cm}       
    \caption{Overview over the methods discussed in this article.}
    \label{table:overview}
\end{table}

\section{Unsupervised Models}\label{sec:unsupervised}
We start this survey by considering unsupervised methods, that is, we assume access to an independent and identically distributed (iid) data set $(u_i)_{i=1}^N$ of images following the prior distribution $\p(u)$. The first two subsections \Cref{sec:learned_generators,sec:learned_prior} elaborate on approaches for training an approximate model of $\p(u)$, respectively $\Rc(u)$ such that $\p(u) \propto \exp(-\Rc(u))$. After we have access to such a model, the ground truth $u$ for a new observation $y$ can be estimated either by minimizing \eqref{eq:posterior} or by sampling from the posterior $\p(y|u)$ as explained in the introduction. Then, in \Cref{sec:pnp,sec:red}, we will discuss Plug and Play approaches and regularization by denoising. In the former, generic denoisers are used in place of the proximal mapping ${\prox}_{\Rc}(u):= \argmin _v \frac{\|u-v\|_2 ^2}{2} + \Rc(v)$ appearing in optimization algorithms, and in the latter denoisers replace the gradient $\nabla\Rc$, such that MAP estimation or inference based on $\p(u|y)$ is possible without knowing $\Rc$ explicitly.

\subsection{Learned Generators}\label{sec:learned_generators}
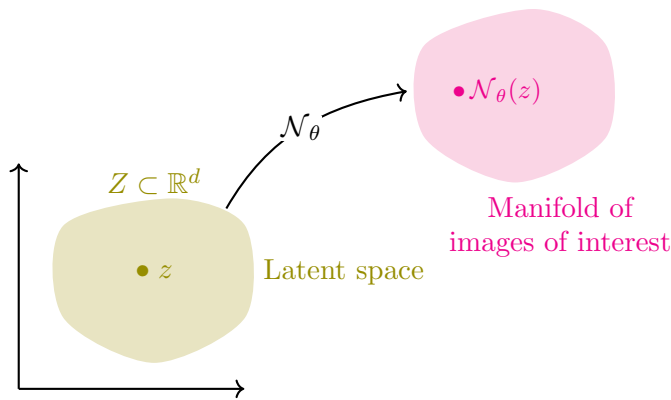
\begin{figure}
\centering
\begin{tikzpicture}[scale = 1.2]
\begin{scope}[every node/.style={inner sep=0,outer sep=0}]
	\draw[magenta!20,thick,fill=magenta!20,xshift=1cm] plot[smooth cycle]
  coordinates {(0+4,1+2) (1+4,0.5+2) (2+4,1+2) (2+4,2+2) (1.2+4,2.4+2) (0+4,2+2)};

   \draw[olive!20,thick,fill=olive!20,xshift=1cm] plot[smooth cycle]
  coordinates {(0,1) (1,0.5) (2,1) (2,2) (1.2,2.3) (0,2)};

  \draw[thick,->] (.5,0.2) -- (3,0.2); 
  \draw[thick,->] (.5,0.2) -- (.5,2.7); 
  \node[olive] (U) at (2,2.5) {$Z\subset\mathbb{R}^d$};
  \node[olive] (mu) at (2,1.5) {$\bullet \;z$};
  \draw [->, thick,out=60,in=190] (2.8,2.2) to node[midway,fill=white]{$\Nc_\theta$} (4.8,3.5);
  \node[olive] (l) at (4.1,1.5) {Latent space};
  
    \node[magenta,align=center] (m) at (6.5,2) {Manifold of \\ images of interest};%
    \node[magenta] (NNmu) at (5.8,3.5) {\small $\bullet \Nc_\theta(z)$};

	\end{scope}
\end{tikzpicture}
\caption{Illustration of a learned generator.}
\label{fig:learned_generator}
\end{figure}

In this article \emph{learned generators} refers to models where one aims to learn a parametrization of the manifold of  images of interest, cf. \Cref{fig:learned_generator}. More precisely, given a tractable latent distribution $\p(z)$, e.g., a Gaussian distribution, a model $\Nc_\theta$ is trained such that for $z\sim\p(z)$ it holds $\Nc_\theta(z)\sim \p(u)$. That is, the network pushes the latent distribution forward to the image distribution. In that sense, $\p(u)$ is not learned explicitly, but rather implicitly via the parametrization $\Nc_\theta$. Employing this parameterization, the inverse problem \eqref{eq:inverse_problem} is transferred to
\begin{equation}
y = \fwd \Nc_\theta(z) + \zeta,
\end{equation}
and we find for the posterior distribution of the latent variable
\begin{equation}
-\log\p(z|y) = -\log \p(y|z) - \log \p(z) + const.
\end{equation}
In particular, in the case of Gaussian noise $\zeta\sim\Nc(0,\sigma^2)$ and a Gaussian distribution on the latent space $z\sim\Nc(0,1)$ we obtain
\begin{equation}\label{eq:learned_generator_posterior}
-\log\p(z|y) = \frac{1}{2\sigma^2} \|\fwd\Nc_\theta(z) - y\|_2^2 + \frac{1}{2}\|z\|_2^2 + const.
\end{equation}
Consequently, \eqref{eq:learned_generator_posterior} allows to apply a learned generator in a variational as well as a Bayesian setting. A crucial degree of freedom in such methods utilizing a learned generator is the approach for training $\Nc_\theta$, for which a maximum likelihood estimation, generative adversarial networks (GAN), or (variational) autoencoders (VAE) are applicable.

\subsubsection{Maximum Likelihood Estimation}
In order to perform maximum likelihood estimation, we first need to obtain an explicit formulation for the likelihood induced by the generator. This, however, can only be achieved if the generator $\Nc_\theta$ is a diffeomorphism, that is, continuously differentiable and invertible with continuously differentiable inverse \cite{kingma2018glow,dinh2017density,dinh2015nice,karami2019invertible}. This poses extensive constraints on the architecture, in particular, prohibiting the use of ReLU activations and imposing restrictions on linear/convolutional layers of the network. Performing a change of variables, we find for $\p_\theta$ the distribution of $\Nc_\theta(z)$
\begin{equation}
\p_\theta(u) = \p(z) \left| D \Nc_\theta^{-1}(u) \right|
\end{equation}
where $z=\Nc_\theta^{-1}(u)$ and $D \Nc_\theta^{-1}$ denotes the Jacobian of $\Nc_\theta^{-1}$. Such models, transforming a complicated distribution to a tractable one via a composition of invertible mappings, e.g., an invertible neural network, are referred to as \emph{normalizing flows} in the literature \cite{kingma2018glow,rezende2015variational}. Given an iid sample of images following the prior distribution $(u_i)_{i=1}^N$, $u_i\sim\p(u)$, the explicit expression for $\p_\theta$ allows to pursue maximum likelihood estimation of the parameters $\theta$ by solving
\begin{equation}
\argmax\limits_{\theta} \sum\limits_{i=1}^N \log \p_\theta(u_i).
\end{equation}
An application of normalizing flows for inverse problems as explained can be found in \cite{asim2020invertible}. Another related approach to employ normalizing flows to inverse problems is to train the flow directly as a parametrization of the posterior $\p(u|y)$ rather than the prior \cite{sun2021deep,rizzuti2020param,whang2021composing}. Note, however, that these methods assume an already given prior distribution $\p(u)$. Thus, one either has to rely on a handcrafted or a separate pre-trained prior. On the other hand, it is possible to learn a \emph{conditional normalizing flow}, that is, a network with inputs $z$ \emph{and} $y$ such that for a fixed $y$ it holds that $\Nc_\theta(y,z)\sim\p(u|y)$ for $z\sim\p(z)$. Training such a conditional flow, however, typically leads to a supervised approach \cite{winkler2019learning,siahkoohi2020faster}.

\subsubsection{Autoencoders}
In the case of autoencoders \cite{obmann2020deep}, we aim to find two transformations, an encoder $E_\phi:\Xc\rightarrow\Zc$ which maps images to latent variables and a decoder $D_\theta:\Zc\rightarrow\Xc$ which, in reverse, maps latent variables to images. Subsequently, the decoder is used as the learned generator. In practice, usually the latent space is chosen low-dimensional compared to the image space, $\text{dim}(\Zc)<<\text{dim}(\Xc)$, which acts as a regularization \cite{obmann2020deep}. We wish to find parameters $(\phi,\theta)$ of encoder and decoder such that $D_\theta\circ E_\phi(u)\approx u$ for any image from our distribution $\p(u)$. Thus, the corresponding training problem reads, e.g., as
\begin{equation}\label{eq:autoencoder_training}
\argmin\limits_{\theta,\phi} \frac{1}{N}\sum\limits_{i=1}^N \ell(D_\theta(E_\phi(u_i)),u_i)
\end{equation}
where $\ell(\;\cdot\;,\;\cdot\;)$ denotes some loss functional, e.g., the MSE loss $\ell(u,v) = \frac{1}{2}\|u-v\|_2^2$. Of course adding regularization of the parameters $(\phi,\theta)$, e.g., via penalizing their norms in \eqref{eq:autoencoder_training} is possible and reasonable to ensure well-posedness of the training problem.

\subsubsection{Variational Autoencoders}
While in practice variational autoencoders (VAE) \cite{kingma2013auto} are implemented using a decoder and an encoder network in a similar fashion to \emph{non-variational} autoencoders, the underlying mathematical formulations of the two differ significantly. Again, our goal is to determine the parameters $\theta$ of a complex distribution $\p_\theta(u)$ such that $\p_\theta(u)\approx \p(u)$. We assume that $u$ depends on the latent variable $z$. The latent distribution $\p(z)$ of $z$ is modeled as a standard Gaussian. The conditional distribution $\p_\theta(u|z)$, on the other hand, is chosen as a multivariate Gaussian with mean and covariance $\mu_z$ and $\Sigma_z$, which are functions of $z$, both modeled as neural networks with parameters $\theta$. That is, in this case the neural networks map the latent variable $z$ to the parameters of a distribution and not to a data point. The variance $\Sigma_z$ is frequently restricted to be a diagonal matrix which reduces complexity significantly \cite[Appendix C.2]{kingma2013auto}. In subsequent computations we will also make use of the distribution of $z|u$. However, since this distribution is not tractable in general based on the model $\p_\theta(u|z)$ we approximate it with a second model $\q_\phi(z|u)$. As for $\p_\theta(u|z)$, $\q_\phi(z|u)$ is modeled as a Gaussian with parameters $\mu_u$ and $\Sigma_u$, which are neural network functions with parameters $\phi$. For training we aim to maximize the log-likelihood (or \emph{evidence}) $\sum_i\log \p_\theta(u_i)$ of the model with respect to $\theta$. However, this approach is not directly tractable which is why we employ a lower bound of the likelihood referred to as the evidence lower bound (ELBO) derived as follows:
\begin{equation}
\begin{aligned}
\log \p_\theta(u) = \log \int \p_\theta(u,z) \;\d z &= \log \int \frac{\p_\theta(u,z)}{\q_\phi(z|u)} \q_\phi(z|u) \;\d z \\
&\geq \int \log \left[ \frac{\p_\theta(u,z)}{\q_\phi(z|u)} \right] \q_\phi(z|u)\;\d z\\
&= \E_{z\sim \q_\phi(z|u)}\left[ \log \p_\theta(u,z) - \log \q_\phi(z|u)\right]\\
&= \E_{z\sim \q_\phi(z|u)}\left[ \log \p_\theta(u|z) + \log \p(z) - \log \q_\phi(z|u)\right]\\
&= \E_{z\sim \q_\phi(z|u)}\left[ \log \p_\theta(u|z) \right] -\KL\left(\q_\phi(z|u)\middle|\middle|\p(z)\right)\\
&\eqqcolon \text{ELBO}(\theta,\phi,u)
\end{aligned}
\end{equation}
where we used Jensen's inequality in the second line and $\KL\left(\q\middle|\middle|\p\right):= \int \log(\frac{\q(z)}{\p(z)}) \q(z)\d z$ is the Kullback-Leibler (KL) divergence with $\frac{\q(z)}{\p(z)}$  the Radon Nikod\'ym derivative for the distributions $\p,\q$ such that $\q$ is absolutely continuous w.r.t $\p$. Given $(\theta,\phi,u)$, computation of $\text{ELBO}(\theta,\phi,u)$ is feasible and simple since the KL divergence between Gaussians is explicit and sampling from $z\sim \q_\phi(z|u)$ is straightforward as well given the Gaussian model. Maximizing the likelihood of the entire iid sample, the parameters are then determined as
\[\argmax\limits_{\theta,\phi}\sum\limits_{i=1}^N\text{ELBO}(\theta,\phi,u_i)\]
and we refer to \cite{kingma2013auto} for details on the optimization procedure. Note that given our model of $\p_\theta$ the term $\log \p_\theta(u|z)$ within the ELBO reads as
\begin{equation}\label{eq:vae}
\log \p_\theta(u|z) = -\frac{\|u-\mu_z\|_2^2}{2\sigma_z^2}+const.
\end{equation}
Let us denote the neural network modeling $\mu_z$ as $\mu_z = \Nc_\theta(z)$. Thus, while conceptually $\Nc_\theta$ outputs parameters of a distribution, the fact that during training we maximize \eqref{eq:vae} shows that $\Nc_\theta$ is trained to reconstruct samples from the target distribution $\p(u)$ which motivates using it as a generator for data points as well. In the following we elaborate on works using VAEs as learned generators.

In \cite{bora2017compressed} Bora et al. apply a generator obtained as the decoder of a VAE (or as the generator of a GAN, see below \Cref{sec:gans}) to compressed sensing by minimizing the posterior \eqref{eq:learned_generator_posterior}. They additionally show that under certain conditions the obtained results are close to the ground truth with high probability if the forward operator of the inverse problem is a random Gaussian matrix. Similarly, in \cite{hand2018phase} a generator from a VAE (and a GAN) is used for the non-linear inverse problem of phase retrieval.

In \cite{asim2020blind,asim2018solving} the authors use generators obtained from VAEs (or GANs, \Cref{sec:gans}) for the solution of blind image deblurring. They use two distinct generators, one for the blur kernel and one for the sought image. In \cite{dhar2018modeling} the authors follow a similar approach, however they allow for results with sparse deviations from the range of the generator.

In \cite{gonzales2022solving,duff2023vaes} the authors consider MAP estimation for the posterior of ground truth $u$ \emph{and} latent variable $z$, $\p_\theta(u,z|y) \propto \p(y|u,z)\p_\theta(u|z)\p(z)$. That is, contrary to the other works, they do not use the VAE as a deterministic mapping $z\mapsto u$ but rather take into account that according to the VAE formulation the generator merely yields distribution parameters. Moreover, in \cite{duff2023vaes} the authors model $\Sigma_z^{-1} = L_zL_z^T$ with $L_z$ a lower triangular matrix with positive diagonal entries. Additionally, sparsity constraints are imposed on $L_z$.

In \cite{seo2019learning} it is proposed to solve the inverse problem of electrical impedance tomography with a learned approach in two steps. The authors train an encoder $E_{\phi_1}$ and a decoder $D_{\phi_2}$ in a VAE scheme on samples $(u_i)_i$. Afterwards, in a supervised fashion, an additional NN $\Nc_\theta$ is trained to map the the observation to the encoded variables, i.e., $\Nc_\theta(y_i)\approx E_{\phi_1}(u_i)$. After training, for a new observation $y$, an estimate of $u$ is then obtained as $u=D_{\phi_2}(\Nc_\theta(y_i))$. That is, the prediction method itself is a supervised machine learning method, however, a VAE trained in an unsupervised fashion is used to transfer the problem to latent space (and back).

\subsubsection{Generative Adversarial Networks}\label{sec:gans}
In the case of generative adversarial networks (GAN) \cite{goodfellow2014generative} we employ two neural networks during training, a generator $G_\theta$ and a discriminator $D_\phi$ (not to be confused with a decoder) which pursue opposing goals. The generator aims to generate images of the desired target distribution $\p(u)$. The discriminator, on the other hand, tries to distinguish samples from the distribution $\p(u)$ from samples generated by the $G_\theta$. More precisely, $D_\phi(u)$ yields a number in $[0,1]$ which represents the estimated probability of its input $u$ being a \emph{real} sample, that is, the discriminator's goal is to predict $D_\phi(u)\approx 1$ for $u\sim\p(u)$ and $D_\phi(G_\theta(z))\approx 0$, whereas the generators goal is to achieve $D_\phi(G_\theta(z))\approx 1$. A mathematical formulation of the training problem reads, e.g., as
\begin{equation}\label{eq:gan_training}
\min\limits_{\theta}\max\limits_{\phi} \E_{z\sim\p(z)}\left[\log \left(1-D_\phi(G_\theta(z))\right)\right] + \E_{u\sim\p(u)}\left[\log\left(D_\phi(u)\right)\right].
\end{equation}
As mentioned above, the works \cite{bora2017compressed,asim2020blind,dhar2018modeling,asim2018solving,hand2018phase} perform their experiments with generators obtained from VAEs as well as from GANs. Specifically, they use the deep convolutional GAN (DCGAN) from \cite{radford2016dcgan}.

In \cite{shah2018solving}, instead of optimization over the latent space by minimizing $\|\fwd G_\theta(z)-y\|_2^2$ similar to Bora et al. \cite{bora2017compressed}, the authors consider solving the constraint minimization problem
\[\argmin\limits_u \|\fwd u-y\|_2^2,\quad \text{s.t.} \;u\in \rg (G_\theta)\]
using projected gradient descent. Here $\rg (G_\theta)=\left\{G_\theta(z)\;\middle|\;z\in\R^n\right\}$ denotes the range of the generator. In \cite{raj2019ganbased} this approach is further modified by replacing the inner loop necessary for the projection onto the range of the generator by a learned projection which is parametrized as an additional neural network.

In \cite{pan2022exploiting} the authors use a pre-trained GAN generator and propose to solve 
\begin{equation}
\min_{z,\theta} \Dc_y(\fwd G_\theta(z)),
\end{equation}
that is, they also allow for additional fine tuning of the parameters of the generator during inference similar to the deep image prior \cite{ulyanov2018deep}, see \Cref{sec:untrained}. Moreover, as a data fidelity term they use 
\begin{equation}
\Dc_y(u) = \sum\limits_k \|D_\phi(y,k)-D_\phi(\fwd u,k)\|_1
\end{equation}
where $D_\phi(y,k)$ denotes the $k$-th layer of the discriminator network.

In \cite{yeh2017semantic} the problem of large hole inpainting is tackled using a GAN generator. Instead of the frequently used squared $L^2$ loss $\frac{1}{2}\|z\|_2^2$ which stems from thestandard Gaussian latent distribution model $\log\p(z)=-\frac{1}{2}\|z\|_2^2+const.$, the authors propose using $\log (1-D_\phi(G_\theta(z)))$ as a penalty on the latent variable $z$ during MAP estimation.

\noindent\textbf{Wasserstein GAN:}
It can be shown that \eqref{eq:gan_training} amounts to minimizing a specific distance, namely the Jensen-Shannon divergence, between the target density and the the density generated by $G_\theta$ \cite{goodfellow2014generative}. In \cite{arjovsky2017wasserstein} the authors propose to replace the Jensen-Shannon divergence with the Wasserstein-1 distance which improves numerical behavior of the training. Using Kantorovich-Rubinstein duality \cite{villani2009optimal} the resulting training problem reads as
\begin{equation}\label{eq:wasserstein_gan_training_1}
\min\limits_{\theta}\max\limits_{D_\phi\in\text{$1$-Lip}} \E_{u\sim\p(u)}\left[D_\phi(u)\right]-\E_{z\sim\p(z)}\left[ D_\phi(G_\theta(z)) \right]
\end{equation}
where the discriminator $D_\phi$ is constrained to be Lipschitz continuous with Lipschitz constant $1$. This constraint is, however, difficult to enforce strictly, which is why it is commonly relaxed by replacing the discriminator update for a fixed generator with \cite{gulrajani2017improved}
\begin{equation}\label{eq:wasserstein_gan_training}
\max\limits_{\phi} \E_{u\sim\p(u)}\left[D_\phi(u)\right]-\E_{z\sim\p(z)}\left[ D_\phi(G_\theta(z)) \right] - \lambda\E_{\hat{u}\sim\p(\hat{u})}\left[ (\|\nabla_u D_\phi(\hat{u})\|_2-1)^2 \right]
\end{equation}
where $\p(\hat{u})$ denotes a uniform distribution on straight lines between a sample from $\p(u)$ and a sample generated from $G_\theta$. This training strategy, referred to as Wasserstein GAN is, e.g., used in \cite{mosser2020stochastic} for sampling from the posterior using a GAN for the specific problem of seismic waveform inversion. Wasserstein GANs are also frequently applied in supervised settings and we refer to \Cref{sec:post-processing} for more info.

\subsection{Learned Priors}\label{sec:learned_prior}
In contrast to the methods in \Cref{sec:learned_generators}, where the prior is implicitly learned via learning the pushforward from a tractable distribution to $\p(u)$, it is also possible to directly approximate the prior, respectively regularizer. An early and prominent example thereof is Fields of Experts (FoE) \cite{roth2005fields,ro09} where the prior is modeled as a product of simple distributions with trainable parameters. In this work, we focus on the approach of parametrizing the regularizer directly as a neural network, i.e., $\Nc_\theta(u) = -\log \p_\theta(u) \approx -\log \p(u)$ \cite{zach2022computed,zach2023stable}, but note that there exist of course many more subtle approaches such as \cite{altekruger2022patchnr}, which trains a patch-based prior using invertible neural networks that is particularly dedicated to limited training data.

A natural approach for directly training a neural network as a regularization functional is the already mentioned maximum likelihood based estimation according to 
\begin{equation}\label{eq:ebm_ML}
\argmax\limits_{\theta} \sum\limits_{i=1}^N \log \p_\theta(u_i)
\end{equation}
as pursued in \cite{zach2022computed,zach2023stable,guan2023magnetic,zongjiang2023kspace} for computed tomography and MRI. 
It is worth noting that, in the limit for $N\rightarrow\infty$ according to the law of large numbers, maximum likelihood estimation as above is, in fact, equivalent to minimizing the Kullback-Leibler divergence between $\p(u)$ and $\p_\theta(u)$ if we scale \eqref{eq:ebm_ML} by $\frac{1}{N}$
\begin{equation}\label{eq:ebm_ML2}
\begin{aligned}
\lim\limits_{N\rightarrow\infty}\frac{1}{N}\sum\limits_{i=1}^N \log \p_\theta(u_i) = &\E_{u\sim\p(u)}[\log \p_\theta(u)]\\
=&\E_{u\sim\p(u)}[\log \p(u)]-\E_{u\sim\p(u)}[\log \frac{\p(u)}{\p_\theta(u)}] \\
= &-\KL\left(\p(u)\middle|\middle|\p_\theta(u)\right)+const.
\end{aligned}
\end{equation}

\subsubsection{Score-based and diffusion models}\label{sec:score-based}
Maximum likelihood estimation according to \eqref{eq:ebm_ML} is commonly performed using gradient based optimization techniques which require the gradient $\nabla_\theta \E_{u\sim\p(u)}[\log \p_\theta(u)]$. Identifying once again the prior with a learned regularizer via
\begin{equation}
\p_\theta(u) =\frac{\exp(-\Rc_\theta(u))}{\int\exp(-\Rc_\theta(v))\;\d v}
\end{equation}
we can compute the gradient using Leibniz's rule for differentiating an integral \cite{zach2023stable,hinton2002training}
\begin{equation}\label{eq:ebm_gradient}
\nabla_\theta \E_{u\sim\p(u)}[\log \p_\theta(u)] = -\E_{u\sim\p(u)}[\nabla_\theta\Rc_\theta(u)] + \E_{u\sim\p_\theta(u)}[\nabla_\theta \Rc_\theta(u)].
\end{equation}
Unfortunately, this reveals a downside of this approach, namely that the second term in \eqref{eq:ebm_gradient} is an expected value with respect to $\p_\theta(u)$ which is computationally demanding as it has to be approximated using sampling procedures.

\noindent\textbf{Score-based modeling:}
A remedy is provided by score-based methods: Interestingly, whether we try to compute the MAP via minimizing the variational objective functional or the MMSE estimate via sampling from the posterior, in any case we do not require access to the value of the regularizer, respectively prior, $\log \p(u)$ but only its gradient $\nabla\log\p(u)$ called the \emph{score}. This gives rise to a class of methods, referred to as \emph{score-based} modeling, where one aims to directly approximate $\nabla\log\p(u)$ using a neural network $s_\theta(u)\approx \nabla\log\p(u)$ \cite{hyvarinen2005estimation,song2019generative,song2020sliced}. The corresponding training problem reads as
\begin{equation}\label{eq:score_matching}
\min\limits_\theta \J(\theta)\coloneqq \E_{u\sim\p(u)}\left[\frac{1}{2}\|s_\theta(u)-\nabla \log\p(u)\|_2^2\right].
\end{equation}
While this functional is not amenable to optimization as we do not have access to $\nabla \log\p(u)$, via integration by parts it can be shown \cite{vincent2011connection} that 
\begin{equation}\label{eq:score_matching2}
\J(\theta)= \E_{u\sim\p(u)}\left[\frac{1}{2}\|s_\theta(u)\|_2^2 + \tr \left(\nabla_u s_\theta(u)\right)\right] + const.
\end{equation}
which can be minimized.

\noindent\textbf{Denoising score matching:}
In \cite{song2019generative}, however, the authors argue that direct estimation of the score using samples from $\p(u)$ unfortunately suffers from several issues in the case of imaging: It is reasonable to assume, that $\p(u)$ is supported in a lower-dimensional manifold rendering the gradient $\nabla \log \p(u)$ unamenable. Moreover the estimation of the score will be poor in areas of low probability and lastly, multi-modality of $\p(u)$ renders sampling difficult. As a remedy it is proposed to consider smoothing the distribution $\p(u)$ with additive Gaussian noise $\tilde{u} = u+\sigma_i\zeta$, $\zeta\sim\Nc(0, \I)$, $0<\sigma_L<\dots<\sigma_1$ leading to the smoothed distributions $\p_{\sigma_i}(\tilde{u}) = \p*\Nc(0,\sigma_i^2 \I)(\tilde{u})$, that is, convolutions of $\p(u)$ with Gaussians of different variances. For these distributions all the aforementioned issues are reduced. The \emph{noise conditional} score network $s_\theta(\tilde{u},\sigma_i)$ is then trained to approximate the score $\nabla\log\p_{\sigma_i}(\tilde{u})$ for all $\sigma_i$, where for $\sigma_i\rightarrow 0$ the desired score should be recovered, by solving
\begin{equation}\label{eq:explicit_score_matching}
\min\limits_{\theta}\sum\limits_{i=1}^L \lambda(\sigma_i) \E_{\tilde{u}\sim\p_{\sigma_i}(\tilde{u})}\left[ \left\| s_\theta(\tilde{u},\sigma_i) - \nabla \log \p_{\sigma_i}(\tilde{u})\right\|_2^2 \right],
\end{equation}
where $\lambda$ is a coefficient function allowing to weigh different noise levels differently. As before, due to $\nabla \log \p_\sigma(\tilde{u})$ being unknown, the formulation \eqref{eq:explicit_score_matching} is unfortunately not amenable to optimization as is. In this case, however, as alternative to \eqref{eq:score_matching2}, the following equivalent problem can be derived \cite{vincent2011connection}:
\begin{equation}\label{eq:denoising_score_matching}
\min\limits_{\theta}\sum\limits_{i=1}^L \lambda(\sigma_i) \E_{u\sim\p(u)}\E_{\tilde{u}\sim\p_{\sigma_i}(\tilde{u}|u)}\left[ \left\| s_\theta(\tilde{u},\sigma_i) - \nabla \log_{\tilde{u}} \p_{\sigma_i}(\tilde{u}|u)\right\|_2^2 \right].
\end{equation}
That is, we can replace the score $\nabla \log \p_{\sigma_i}(\tilde{u})$ by the conditional one $\nabla_{\tilde{u}} \log \p_{\sigma_i}(\tilde{u}|u)$ which is explicit. %
Note that by definition of the smoothed distributions, $\nabla_{\tilde{u}} \log \p_\sigma(\tilde{u}|u) = \frac{1}{\sigma^2}(u-\tilde{u})$, and, thus, 
\[\sigma^2\left( s_\theta(\tilde{u},\sigma) - \nabla_{\tilde{u}} \log \p_\sigma(\tilde{u}|u)\right)=(\tilde{u}+\sigma^2s_\theta(\tilde{u},\sigma))-u\]
showing that $s_\theta(\tilde{u},\sigma)$ is trained to satisfy $(\tilde{u}+\sigma^2s_\theta(\tilde{u},\sigma))\approx u$, thus, resembling a denoiser.

After the score is trained, sampling using the noise conditional score network is performed using the so-called \emph{annealed Langevin} algorithm \cite{song2019generative} which refers to the procedure of consecutively sampling from $\p_{\sigma_i}$ for decreasing noise levels $i=1,\dots,L$ and for each new $i$ using the last iterate of the previous Langevin sampling as initialization. That is, we generate Markov chains according to
\begin{equation}
u_i^{k+1} = u_i^k + \gamma_i s_\theta(u_i^k,\sigma_i) + \sqrt{2\gamma_i}\zeta^k_i,\quad k=0,\dots, K-1
\end{equation}
where $\zeta^k_i$ denote iid standard Gaussian random variables and for $i=1$ the chain is initialized arbitrarily and for $i>1$ we use $u_i^0=u_{i-1}^K$. The step sizes $\gamma_i$ constitute hyperparameters of the method and can, for instance, be chosen to decrease as the noise scale decreases \cite{song2019generative}.

\noindent\textbf{Diffusion models:}
A crucial ingredient in score-based modeling as introduced above is to \emph{corrupt} the target distribution with noise leading to a more well-behaved distribution $\p_{\sigma_i}$. Afterwards sampling from $\p(u)$ is performed by consecutively sampling from $\p_{\sigma_i}$, $i=1,\dots,L$. In other words, \emph{we start at a noisy distribution and successively remove the noise from the samples until the target distribution is obtained}. This approach is amenable for generalization. Let us define the noise corruption process starting at the target distribution $\p(u)$ and leading to the well-behaved distribution as the following Markov chain
\begin{equation}\label{eq:diffusion_discrete}
u_{i+1} = \alpha_i u_i + \beta_i\zeta_i,\quad i=0,\dots,L-1
\end{equation}
with $u_0\sim\p(u)$. By choosing $\alpha_i=1$, $\beta_i=\sqrt{\sigma_{i+1}^2-\sigma_i^2}$, and $\sigma_0=0$ we recover exactly the distributions $u_i\sim\p_{\sigma_i}(u_i)$ previously introduced in the context of denoising score matching. This choice is referred to as \emph{variance-exploding} (VE) diffusion. Another popular choice is to pick $\beta_i\in (0,1)$ and $\alpha_i=\sqrt{1-\beta_i^2}$ such that the variance $\var [u_i] = 1$ for all $i$ if $\var [u_0]=1$. Naturally, this approach is referred to as \emph{variance-preserving} (VP).

As shown in \cite{song2020score} the Markov chain \eqref{eq:diffusion_discrete} admits a time-continuous formulation as a stochastic differential equation (SDE) of the form
\begin{equation}\label{eq:sde}
du_t = f(t)\d t + g(t)\d B_t,\quad 0\leq t\leq T
\end{equation}
with the initial condition $u_0\sim\p(u)$. Here $B_t$ denotes Brownian motion. We denote the distribution at time $t$ as $u_t\sim\p_t(u)$. As $L\rightarrow\infty$ in \eqref{eq:diffusion_discrete} the VE diffusion scheme leads to the time-continuous formulation \cite[Appendix B]{song2020score}
\begin{equation}
f(t)\equiv 0, \quad g(t) = \sqrt{\frac{\d \sigma^2(t)}{\d t}}
\end{equation} 
with an increasing function $\sigma(t)$ such that the final distribution satisfies $\p_T(u)\approx \Nc(0,\sigma(T)^2\I)$. The VP diffusion scheme, on the other hand, results in 
\begin{equation}
f(t)= -\frac{\beta(t)}{2}, \quad g(t) = \sqrt{\beta(t)}
\end{equation} 
and $\p_T(u)\approx \Nc(0,\I)$. Recall that our goal eventually is to revert the process \eqref{eq:sde}. To do so we make use of a result shown in \cite{anderson1982reverse} stating that the SDE \eqref{eq:sde} admits a \emph{reverse SDE}, that is, an SDE which flows backwards in time and leads to a process of the same distribution as \eqref{eq:sde}. The reverse SDE reads as 
\begin{equation}\label{eq:sde_backward}
du_t = \left[f(t)-g(t)^2\nabla \log \p_t(u)\right]\d t + g(t)\d \bar{B}_t
\end{equation}
where now $\d t$ is a negative infinitesimal time step and $\bar{B}_t$ denotes Brownian motion moving backwards in time.

A natural way of sampling from $\p(u)$ using this framework is to discretize and simulate \eqref{eq:sde_backward} from $t=T$ to $t=0$ starting from the tractable distribution $\p_T(u)$ \cite{song2020score,song2021solving,ho2020denoising}. In order to do so, in a training stage we first have to estimate the score $\nabla \log \p_t(u)$ which is possible using denoising score matching \eqref{eq:denoising_score_matching}.

Interestingly, it is also feasible to pursue a slightly different modeling approach which, however, turns out equivalent in the end. Note that \eqref{eq:sde_backward} shows that the reverse diffusion is of the same functional form as the forward one. Thus, if $\p(u_{i+1}|u_i)=\Nc(u_{i+1};\alpha_i u_i,\beta_i^2 I)$ denotes the transition probabilities of the forward diffusion according to \eqref{eq:diffusion_discrete}, the transition probabilities of the reverse diffusion, denoted as $\q(u_i|u_{i+1})$, will also be Gaussian $\q(u_i|u_{i+1})=\Nc(u_i;\mu(u_{i+1},i+1),\Sigma(u_{i+1},i+1))$ where $\mu(u_{i+1},i+1),\Sigma(u_{i+1},i+1)$ are functions of the state $u_{i+1}$ and the time step $i+1$. In practice the variance might be assumed to be fixed fixed, e.g., as $\Sigma(u_{i+1},i+1)=\beta_{i+1}$ \cite{ho2020denoising}. Thus, instead of training the score, we can instead as well train the means of the reverse diffusion $\mu_\theta(u_{i+1},i+1)\approx\mu(u_{i+1},i+1)$ using a neural network $\mu_\theta$. To be precise, the approximations $\mu_\theta(u_{i+1},i+1)$ lead to transition probabilities $\q_\theta(u_i|u_{i+1})$ and eventually to the distribution $\q_\theta(u_0)$ which we can use for maximum likelihood estimation, that is, maximizing $\theta\mapsto\E_{u\sim\p(u)}\left[ \log\q_\theta(u_0) \right]$ so that $\q_\theta(u_0)\approx\p(u)$. A lower bound $L(\theta)$ for this objective is derived in \cite{sohl2015deep,luo2023bayesian} as
\begin{equation}\label{eq:diffusion_MLE}
L(\theta) = \sum_{i>1}\E_{u_0\sim\p(u)}\left[ \KL\left(\p(u_{i-1}|u_i,u_0)\|\q_\theta(u_{i-1}|u_i)\right) \right]+const.
\end{equation}
It turns out that, in the case of Gaussian transitions $\p(u_{i+1}|u_i)$, training via \eqref{eq:diffusion_MLE} is, indeed, equivalent to training a score $s_\theta(u_i,i)\approx \nabla\log \p_{t_i}(u_i)$ via denoising score matching. The approximated score can be used to compute the means $\mu_\theta(u_{i+1},i+1)$ of the reverse diffusion transition probabilities and vice versa \cite{ho2020denoising,luo2023bayesian}.

An interesting work worth mentioning in this context is also \cite{zach2023explicit}, where the authors do not approximate the unknown density $\p_t(u)$ using a general neural network, but show analytically that certain Gaussian mixture models solve the diffusion process exactly. The remaining parameters of these models are then trained using denoising score matching.

Having access to the trained score $s_\theta(u_i,i)\approx \nabla\log \p_{t_i}(u_i)$, some authors propose to complement the approximation of the reverse diffusion with iterations of a MCMC sampling method for the distribution $\p_{t_i}(u_i)$ for fixed $t_i$. Thus, the output of one step of the reverse diffusion process acts as a proposal which is then used as initialization for multiple steps of an  MCMC method. This technique of alternating (reverse) time steps with MCMC steps for fixed time is referred to as \emph{predictor/corrector sampling} \cite{chung2022score,song2020score}. To reiterate, altogether three prevalent sampling methods appear within the realm of diffusion based methods: 1) Annealed Langevin sampling \cite{song2019generative}, i.e., MCMC sampling from successively less noisy distributions, 2) purely discretizing the reverse SDE \cite{song2020score,song2021solving,ho2020denoising}, and 3) as a mixture of the two predictor/corrector sampling \cite{chung2022score,song2020score}

\noindent\textbf{Application to inverse problems:}
By now, we have explained how score-based and diffusion models are used to approximately sample from a complex distribution $\p(u)$. In view of inverse problems we are, however, interested in samples from $\p(u|y)$ where $y$ denotes the given observation. Various techniques to incorporate the observation have been proposed in the literature.

In \cite{chung2022score,erlacher2023joint} the authors apply a diffusion model to MRI and incorporate the observation by performing predictor/corrector sampling for $\p(u)$ and gradient descent steps for the data fidelity in an alternating fashion. In \cite{song2021solving} a stochastic process $y_t = \fwd u_t + \zeta$ for the observation corresponding to $u_t$ is defined and for sampling the authors alternate steps of a discretization of the reverse diffusion with proximal data consistency steps, which balance data consistency to $y_t$ and deviations from the previous iterate. Also in \cite{chung2022come} the strategy of alternating reverse diffusion steps with data consistency steps is employed. The authors additionally propose to initialize the reverse diffusion with the output of the forward diffusion process starting from an initial guess of the solution which can be computed efficiently. This process acts as an acceleration.

In \cite{luo2023bayesian,sohl2015deep} the authors derive expressions for the transition probabilities conditioned on the observation $\q_\theta(u_{i-1}|u_i,y)$ using Bayes theorem. Similarly, in \cite{song2020score}, the authors derive an approximation for the conditional score $\nabla\log \q_\theta(u_t|y)$ given an observation for the inverse problem.

In \cite{chung2022diffusion,chung2022improving} the conditional score is factorized as $\nabla_{u_t} \log \p(u_t|y)=\nabla \log \p(y|u_t)+\nabla_{u_t} \log \p(u_t)$. The latter is estimated using score matching whereas the former is approximated by $\nabla_{u_t} \log \p(y|\hat{u}_0)$ with $\hat{u}_0 = \E[u_0|u_t]$ using Tweedie's formula \cite{efron2011tweedie}. Sampling can then be performed by approximating the reverse diffusion process with the techniques explained above. In \cite{chung2023parallel} this approach is extended to blind inverse problems where the forward operator depends on unknown parameters as well.

In \cite{kawar2021snips}, the authors consider an approach inspired by denoising score matching \cite{song2019generative}. However, they carefully design the noise additions in a way allowing for an explicit computation of the conditional score $\log\p_\sigma(\tilde{u}|y)$ given the unconditional one $\log\p_\sigma(\tilde{u})$ for which they utilize a SVD of the forward operator.

In \cite{kawar2022denoising} the forward and the backward diffusion are defined conditioned on the observation. However, training such a model is undesirable as it depends on the observation $y$. The authors provide a remedy by presenting a result stating that the training of the conditional model is equivalent to that of the unconditional one if there is no weight sharing of networks at different time steps.

\subsection{Plug and Play Priors}\label{sec:pnp}
Recall that MAP estimation is performed by solving
\begin{equation}\label{eq:map_pnp}
\argmin\limits_u \Dc_y(\fwd u)+\lambda\Rc(u)
\end{equation}
with the interpretation $\p(y|u)\propto\exp(-\Dc_y(\fwd u))$ and $\p(u)\propto\exp(-\lambda\Rc(u))$. Different methods for solving \eqref{eq:map_pnp} numerically have been proposed in the literature, such as proximal-gradient (forward-backward) methods \cite{chen1997convergence,bredies2008linear,daubechies2004ista}
\begin{equation}
u^{k+1} = {\prox}_{\tau\lambda \Rc}\left( u^k - \tau\lambda\fwd^*\nabla \Dc_y(\fwd u^k)\right)
\end{equation}
and its accelerations \cite{beck2009fista}, ADMM optimization \cite{boyd2011admm,venkatakrishnan2013plug}
\begin{equation}
\begin{cases}
u^{k+1} = {\prox}_{\frac{1}{\mu}\Dc_y\circ \fwd}(v^k-w^k)\\
v^{k+1} = {\prox}_{\frac{\lambda}{\mu}\Rc}(u^{k+1}+w^k)\\
w^{k+1} = w^k + u^{k+1}-v^{k+1}
\end{cases}
\end{equation}
or primal-dual methods \cite{chambolle2011pd}
\begin{equation}
\begin{cases}
v^{k+1} = {\prox}_{\sigma\Dc_y^*}(v^k+\sigma \fwd\bar{u}^k)\\
u^{k+1} = {\prox}_{\tau\lambda\Rc}(u^k-\tau\fwd^*v^{n+1})\\
\bar{u}^{n+1}=2u^{k+1}-u^k.
\end{cases}
\end{equation}
Similar to the derivation of score-based models, we can observe that non of these methods utilizes the mapping $\Rc$ explicitly, however, they all at some point make use of the proximal mapping
\begin{equation}
{\prox}_{\lambda \Rc}(u)\coloneqq \argmin\limits_{v}\frac{1}{2}\|u-v\|^2_2 + \lambda\Rc(v).
\end{equation}
Note, however, that ${\prox}_{\lambda \Rc}(u)$ is simply the MAP estimate with respect to the observation $u$ in the case of Gaussian noise and with the prior $\p(u) \propto \exp(-\lambda\Rc(u))$. Consequently, in \cite{venkatakrishnan2013plug} it is proposed to replace ${\prox}_{\lambda \Rc}(u)$ with a pre-trained denoiser $D(u)$ mapping noisy to clean images which is referred to as \emph{plug and play} (PnP). This idea of \cite{venkatakrishnan2013plug} motivated a large body of follow up works employing different variations of PnP. In \cite{venkatakrishnan2013plug,brifman2016turning} the PnP approach was based upon the iterations of ADMM. Subsequent works made use of other convex optimization approaches as mentioned, such as FISTA \cite{kamilov2017plug}, half-quadratic splitting \cite{zhang2017learning,zhang2022plug}, or proximal-gradient iterations \cite{ryu2019plug}. Interestingly, in \cite{meinhardt2017learning} it is shown that the fixed points of PnP derived from proximal-gradient iterations, ADMM, or primal-dual iterations are, in fact, identical.

In \cite{tirer2019image} the authors consider the variational problem
\begin{equation}
\min\limits_{u,\tilde{v}} \lambda\|\tilde{v}-u\|_2^2 + \Rc(u),\quad\text{s.t.}\; \fwd\tilde{v}=y
\end{equation}
for which - inspired by PnP - the proximal step is replaced by a denoiser. 

A lot of the subsequent works proposed to use CNN based denoisers within PnP approaches. For instance, in \cite{zhang2017learning,zhang2022plug,ryu2019plug,
tirer2019image,meinhardt2017learning} the authors use CNNs trained for end-to-end denoising either directly or in a residual way such that the CNN reconstructs the noise. A specifically delicate training procedure involving adversarial training of the denoiser is employed in \cite{chang2017one}. In \cite{sun2019online}, on the other hand, the trainable non-linear reaction diffusion (TNRD) denoiser from \cite{chen2017trainable} is used. We want to emphasize at this point that the CNN denoisers used in the context of PnP methods are typically trained to remove noise using pairs of clean and noisy images rendering them as supervised methods conceptually. However, since in the context of PnP methods a pre-trained denoiser can be used for any inverse problem and does not rely on labeled training for the specific problem at hand, we still categorize them as unsupervised in our survey.

In \cite{ye2018deep,dong2019denoising}, on the contrary, the denoiser is trained on image pairs consisting of the noisy filtered back-projection and the ground truth data from CT imaging, rendering this approach fully supervised as also the specific forward model is incorporated in the labeled training data.

The empirical success of PnP methods was succeeded by theoretical considerations. Of course, it is easy to see that in the case that the denoiser is in fact the proximal mapping of some appropriate, convex functional $\Rc$, $D(u)={\prox}_{\Rc}(u)$, convergence of PnP schemes is a simple consequence of convergence of the corresponding convex optimization method (ADMM, proximal-gradient, etc.). This can be shown based on a classical result by Moreau \cite{moreau1965proximite}, under the rather strong conditions that the denoiser has a symmetric Jacobian with eigenvalues in $[0,1]$ \cite{sreehari2016plug}.

In \cite{chan2017plug} convergence for ADMM PnP with an increasing penalty parameter within the ADMM scheme is shown assuming a bounded denoiser, in the sense that $\|D_\sigma(x)-x\|_2^2/n\leq c \sigma^2$ for all noise levels $\sigma$. In \cite{sun2019online} ergodic convergence of ADMM PnP is proven if the denoiser is \emph{averaged}, a property which can, e.g., be ensured by employing a specific network type referred to as \emph{proximal neural networks} \cite{hasannasab2020parseval}. These are neural networks which consist of a concatenation of proximal mappings \cite[Lemma 5.3]{hasannasab2020parseval}. In \cite{ryu2019plug} convergence for proximal-gradient PnP and ADMM PnP is proven using fixed point arguments by enforcing Lipschitz continuity on the residual CNN denoiser.

\subsubsection{Plug and Play via Tweedie}

A slightly different approach to PnP which is closely related to the derivation of score based modeling in \cite{song2019generative} is to approximate the true prior distribution $\p(u)$ by a smoothed version thereof, $\p_\sigma(\tilde{u}) = \p(u)*\Nc (0,\sigma^2 I)(\tilde{u})$, that is, the convolution with a Gaussian kernel. Tweedie's formula then \cite{efron2011tweedie} states that 
\begin{equation}
\E[u|\tilde{u}] = \tilde{u} + \sigma^2\nabla\log\p_{\sigma}(\tilde{u}).
\end{equation}
Thus, if the denoiser satisfies $D_{\sigma}(\tilde{u}) \approx \E [u|\tilde{u}]$, we might approximate $\nabla\Rc(u) = \nabla\log\p_{\sigma}(\tilde{u})$ by
\begin{equation}\label{eq:pnp_tweedie2}
\nabla\log\p_{\sigma}(\tilde{u}) \approx \frac{1}{\sigma^2}\left(D_{\sigma}(\tilde{u}) - \tilde{u}\right)
\end{equation}
and use this approximation within gradient based algorithms for variational inference, respectively MAP estimation \cite{bigdeli2017image,arjomand2017deep}, or for sampling from the posterior $\p_\sigma(\tilde{u}|y)$ and subsequent MMSE or variance estimation  \cite{laumont2022bayesian}. While this only provides a method for inference of a smooth approximation of the target distribution, the discrepancy between $\p_{\sigma}(\tilde{u})$ and $\p(u)$ can be made arbitrarily small for $\sigma\rightarrow 0$ \cite{laumont2022bayesian}. Moreover, in \cite[Theorem 1]{alain2014regularized} it is shown, that also the discrepancy between $\frac{1}{\sigma^2}\left(D_{\sigma}(\tilde{u}) - \tilde{u}\right)$ and the gradient $\nabla\log\p(u)$ tends to zero as $\sigma\rightarrow 0$ if the denoiser is trained minimizing the squared error of the reconstruction for Gaussian noise corruption. In \cite{guo2019agem} such a PnP approach is used to simultaneously estimate the noise level an and an MMSE reconstruction of the inverse problem.

\subsection{Regularization by Denoising}\label{sec:red}
In \cite{romano2017red} the authors proposed \emph{regularization by denoising} (RED), an approach sharing similarities with PnP priors. Specifically, the authors introduce a RED regularizer defined as
\begin{equation}
\Rc(u) = u^T\left( u - D(u)\right)
\end{equation}
where $D(u)$ denotes a denoiser. The denoiser is assumed to satisfy two crucial conditions, namely, i) local homogeneity, that is, $D(cu)=cD(u)$ for $|1-c|<\epsilon$ and ii) stability, i.e., the spectral radius satisfies $\max_i |\lambda_i(\nabla D(u))|<1$ for all $u$ where $\lambda_i$ denote the eigenvalues of the Jacobian of the denoiser. These conditions are empirically demonstrated for several denoising methods (K-SVD \cite{elad2006image,mairal2009nonlocal}, BM3D \cite{dabov2007image}, NLM \cite{buades2005nonlocal,kervrann2006optimal,tasdizen2009principal}, EPLL \cite{zoran2011from}, and TNRD \cite{chen2017trainable}) where it should be mentioned that differentiability is not given for EPLL, BM3D, and K-SVD as is.

In \cite{romano2017red} the authors claim that, based on the assumptions, the RED regularizer satisfies $\nabla \Rc(u) = u-D(u)$ which is used in subsequent gradient based approaches for variational inference, such as gradient descent, ADMM, and a fixed point iteration for determining roots of the gradient of the objective functional, that is, $u$ such that
\begin{equation}
0=\nabla_u(\Dc_y(\fwd u) + \lambda\Rc(u)).
\end{equation}
The approach is applied with two different denoisers $D$, a simple median filter and TNRD \cite{chen2017trainable} as a CNN based method.

It was shown later, that the claims made in \cite{romano2017red} are not always true. To be precise, for continuously differentiable denoisers the identity $\nabla \Rc(u) = u-D(u)$ is only satisfied in the case that the Jacobian of the denoiser $D(u)$ is symmetric \cite{reehorst2019red}. Otherwise, there exists no regularizer admitting $u-D(u)$ as its gradient which is a direct consequence of Schwarz's theorem. In \cite{reehorst2019red} it is also empirically verified, that the gradient identity $\nabla \Rc(u) = u-D(u)$ is significantly violated for the TNRD denoiser \cite{chen2017trainable} and the CNN denoiser \cite{zhang2017beyond}. Despite these insights it is shown in \cite{reehorst2019red} that for non-expansive denoisers, the RED fixed point iteration is, in fact, convergent (however, not to a minimizer of a variational regularization problem with a RED regularizer).

Further theoretical analysis was provided in \cite{cohen2021regularization} where the authors consider a slightly different approach. In this work the denoisers are assumed to be demicontractive and satisfy $D(0)=0$. After showing that $D(u)=u$ is equivalent to $\Rc(u) = u^T\left( u - D(u)\right) = 0$ - that is, the fixed points of the denoiser coincide with the roots of the RED regularizer - the authors propose to solve
\begin{equation}
\min\limits_u \Dc_y(u),\quad \text{s.t. } D(u)=u.
\end{equation}
For this optimization problem a provably convergent hybrid steepest descent algorithm is presented. The method (as well as a relaxed version) is tested on several inverse problems with multiple denoisers, in particular, again TNRD \cite{chen2017trainable}.

In \cite{hong2019acceleration} RED is accelerated employing vector extrapolation. While empirically the method performs as expected, achieving the theoretically necessary conditions for convergence "might not be realistic in practice". The authors tested the method TNRD \cite{chen2017trainable}.

In \cite{sun2019block} the RED approach is adapted to a block coordinate structure. Within the proposed gradient descent algorithm different coordinate blocks of $u$ are updated separately chosen at random. The updates within the algorithm admit ergodic convergence to zero in expectation if the denoiser is non-expansive for each coordinate block. In particular, under the stronger assumption that the denoiser is the proximal mapping of a function with bounded subgradient, the objective function value converges to its minimum in expectation. Experimental results are provided with an end-to-end CNN denoiser and other non-CNN based methods (TV and BM3D).

Similarly, in \cite{wu2019online} the authors consider a stochastic gradient descent based algorithm for RED for which they show ergodic convergence to zero of the updates in expectation. In experiments they use a CNN denoiser trained to remove additive Gaussian noise.

An application of RED to phase retrieval using the residual CNN denoiser from \cite{zhang2017beyond} can be found in \cite{metzler2018prdeep}.

\section{Supervised}\label{sec:supervised}
In this section we consider the case of having access to an iid sample $(u_i,y_i)_i$ of ground truth-observation pairs for training. The discussed supervised methods are classified into four categories: i) fully learned, ii) post-processing, iii) unrolling, and iv) learned regularizers. In i) the entire reconstruction $y\mapsto u$ is learned, whereas in ii) we use a trained method to improve upon a model based reconstruction. Methods according to iii) utilize NNs whose architecture resemble well-established iterative methods for solving inverse problems. That is, a layer of the NN might, e.g., share similarities to proximal-gradient iterations for \eqref{eq:var_reg}. In iv) a regularizer is trained similarly as described in \Cref{sec:learned_prior}, however, this time training is performed in a supervised manner.

\subsection{Fully Learned}\label{sec:end-to-end}
The conceptually simplest approach to use neural networks for solving inverse problems might be what we refer to as \emph{fully learned} in this article (cf. \emph{direct Bayesian inversion} \cite[Section 5.1.2]{arridge2019solving}). This terminology refers to training a neural network $\Nc_\theta$ to directly map the observation $y$ to a point estimate of $u$ \cite{thu2018image,zhang2021residual}. The prototypical training problem for fully learned approaches reads as
\begin{equation}\label{eq:training_supervised}
\min\limits_\theta \sum\limits_{i=1}^N \ell(u_i,\Nc_\theta(y_i))
\end{equation}
where $(u_i,y_i)_i$ denotes a given data set of ground truth/observation pairs satisfying $y_i = \fwd u_i + \zeta_i$ with iid noise realizations $\zeta_i$. Of course, \Cref{eq:training_supervised} can be modified to include, e.g., regularization of the parameters $\theta$. In the case of fully supervised approaches, the trained method has to capture not only the prior distribution or the noise but also the forward model $\fwd$. In particular, this means that for any different inverse problem the method has to be trained again, whereas unsupervised approaches can be used for different forward operators as long as the distribution $\p(u)$ is unchanged. The potentially increased complexity due to implicitly learning the (reverse) forward operator, however, is not an issue in the case of image denoising, for which the forward operator reduces to the identity. Thus, denoising is particularly suited for fully learned approaches \cite{jain2008natural}. A lot of works in this context proposed to train residual neural networks. That is, instead of estimating the ground truth $u$, the network estimates the noise, respectively the residual, $\Nc_\theta(y)\approx y-u$ \cite{gurrola2021residual,tian2019enhanced,zhang2017beyond}.

\subsection{Post-processing}\label{sec:post-processing}
In the context of non-trivial forward operators, an approach significantly more popular than fully learned methods is data-driven post-processing. Instead of learning the entire mapping $y\mapsto u$, in this case, the neural network is preceded by a model-based reconstruction. That is, the neural network is trained to remove any remaining artifacts or noise that are still present after the non-learned reconstruction. 

Substantial research has been conducted for post-processing of MRI \cite{quan2018compressed,yang2018dagan,schlemper2018deep,chang2018deep} and CT \cite{jin2017deep,pelt2018improving,chen2017low,kang2017deep} reconstructions where the zero filling solution and the result obtained via filtered backprojection (FBP), respectively, are further improved by a pre-trained NN.

While training of supervised methods is often performed by simply minimizing the reconstruction loss on a given data set \cite{jin2017deep,pelt2018improving,chen2017low}, a lot of works on post-processing propose adversarial training. In this case, well-known GAN training, according to \eqref{eq:gan_training} or \eqref{eq:wasserstein_gan_training}, has to be adapted to fit the post-processing task, e.g., the model parameters are obtained as solutions to
\begin{equation}\label{eq:gan_postprocessing}
\begin{aligned}
\min\limits_{\theta}\max\limits_{\phi} \E_{\tilde{u}\sim\p(\tilde{u})}\left[\log \left(1-D_\phi(G_\theta(\tilde{u}))\right)\right] + \E_{u\sim\p(u)}\left[\log \left(D_\phi(u)\right)\right] \\+ \E_{(u,\tilde{u})\sim\p(u,\tilde{u})}\left[ \ell(G_\theta(\tilde{u})), u) \right]
\end{aligned}
\end{equation}
where $\tilde{u}$ denotes the corrupted images obtained after applying the model based reconstruction to the observation $y$. The first two terms in \eqref{eq:gan_postprocessing} resemble the usual GAN loss. In the third term, $\ell(G_\theta(\tilde{u})), u)$ denotes some kind of discrepancy measure between the reconstruction obtain via post-processing and the actual ground truth. For instance, $\ell$ might be chosen as a simple MSE loss or a perceptual VGG loss \cite{yang2018dagan}. Note that in this formulation the generator input is not a random sample from a pre-specified latent distribution, but the corrupted image data. 

Several variations of GAN training for post-processing have been proposed. In \cite{quan2018compressed} the authors combine the GAN loss with the MSE or MAE loss in the image as well as in the frequency domain for improving upon the zero filling reconstructions for MRI. In \cite{yang2018dagan} additionally a perceptual VGG loss is included. In \cite{yang2018low} a Wasserstein GAN is combined with a VGG loss for post-processing of low-dose CT images. A different approach is pursued in \cite{you2020ct} where the authors propose to simultaneously train two GANs for post-processing of low resolution CT images. One of the generators maps from low to high resolution and the other one in the opposite direction. The authors employ a Wasserstein GAN loss complemented with additional functionals ensuring, in particular, that the two generators are, indeed, inverse functions of each other.

Another viable approach are conditional GANs \cite{isola2017image} where the generator has two inputs, the corrupted data and a random sample $G_\theta(\tilde{u},z)$. Thus, during application of the method, for fixed $\tilde{u}$ sampling \emph{conditioned on $\tilde{u}$} is possible by sampling from $z$.

\subsection{Unrolling}\label{sec:unrolling}
\emph{Unrolling} or \emph{learned iterative methods} (first proposed in \cite{gregor2010learning}) is a way of designing end-to-end networks for the solution of inverse problems in a more deliberate way. Let us consider the output $u^L$ of a generic $L$ layer feed-forward neural network according to
\begin{equation}\label{eq:unrolling_NN}
u^{k+1} = \sigma\left( W_ku^k + b_k \right),\quad k=0,1,2,\dots, L-1
\end{equation}
where $(W_k,b_k)_k$ are the weights and biases of the network, $\sigma$ is the non-linearity, and $u^0$ is its input. Note that, if we choose the weights and biases as $W_k = (I-\tau\lambda\fwd^* \fwd)$, $b_k=\tau\lambda\fwd^*y$ for all $k$, and the non-linearity as $\sigma(u) = {\prox}_{\tau\lambda\Rc}(u)$, \eqref{eq:unrolling_NN} matches exactly $L$ iterations of the proximal gradient algorithm with step size $\tau$
\begin{equation}\label{eq:unrolling_proxgrad}
u^{k+1} = {\prox}_{\tau\lambda \Rc}\left( u^k - \tau\lambda\fwd^*(\fwd u^k - y)\right)
\end{equation}
for solving
\begin{equation}
\min\limits_{u} \frac{1}{2}\|\fwd u - y\|_2^2 + \lambda\Rc(u).
\end{equation}
Thus, we find a close resemblance between the structure of neural networks and iterative methods for variational regularization. This observation is the foundation of unrolling: Instead of choosing an arbitrary network architecture for solving the inverse problem at hand in an end-to-end fashion, it seems more reasonable to make design choices (loosely) imitating iterations of an optimization method. Specifically, we might keep some parts of the iteration \eqref{eq:unrolling_proxgrad} and replace others by trainable building blocks, e.g., in \cite{xiang2021fista} the authors train only the non-linearities. The training approach of unrolled methods is a crucial difference between unrolling and PnP (\Cref{sec:pnp}). While the two frameworks share similarities in the sense that in both cases one aims to replace certain parts of a well-established iterative method with trainable components, unrolled methods are considered end-to-end and trained in a supervised way according to \eqref{eq:training_supervised}. In PnP methods, on the other hand, the learned component is pre-trained in an unsupervised way and simply "plugged" into the iteration. As a consequence, as in general for supervised methods, the approach in unrolling renders training dependent on the forward operator $\fwd$, such that the trained method cannot be applied to a different inverse problem in general.

The principle of unrolling outlined above is of course not restricted to proximal gradient iterations \cite{gregor2010learning,xiang2021fista,wang2016learning}, which were used as an instructive example, but applicable to any optimization method such as gradient descent \cite{adler2017solving,putzky2017recurrent,hammernik2018learning}, ADMM \cite{yang2016deep,yang2020admm,he2019optimizing}, primal-dual algorithms \cite{adler2018learned} or, as in \cite{chen2017trainable}, to the discretization of a reaction-diffusion equation.

A typical modeling question in the context of unrolling is how much of the method is allowed to be trained and which parts are kept as in the underlying iterative method. In \cite{xiang2021fista,zhang2018ista} the gradient step of a proximal-gradient method is kept as is, whereas the proximal mapping is replaced by a trained version. A similar and empirically successful approach is pursued in the learned primal-dual algorithm \cite{adler2018learned}, where the authors replace the proximal mappings of a primal-dual method by trained CNNs. They additionally improve numerical performance by introducing memory to the algorithm, that is, information of previous iterates is included when computing the next update. In \cite{wang2016learning} on the contrary, the authors unroll a proximal gradient method and preserve the proximal steps, whereas they allow for trainable gradient updates.

Another related modeling question is whether to employ \emph{weight sharing}. Weight sharing refers to imposing an additional constraint on the unrolled network ensuring that the weights/non-linearities of different layers are identical \cite{wang2016learning,xiang2021fista}. This technique reduces the number of trainable parameters and is clearly reasonable if the update rule of the underlying iterative method is identical in each iteration as, e.g., in \eqref{eq:unrolling_proxgrad}. An interesting insight in this context is that weight sharing is closely related to recurrent neural networks as having consecutive layers with identical weights might as well be modeled as a single layer with a feedback loop \cite{putzky2017recurrent}. On the other hand, increasing flexibility by allowing different weights at different layers might lead to improved performance \cite{adler2018learned,zhang2018ista,chen2017trainable,hammernik2018learning} at the cost of higher complexity and less interpretability.

\subsection{Learned Priors}
In addition to the above approaches where an image-to-image mapping is learned, there are several approaches for training a regularization functional, respectively prior, $\Nc_\theta(u)=\Rc(u)$ in a supervised way (contrary to the unsupervised approaches in \Cref{sec:learned_prior}). A rather intuitive formulation to achieve this is a bilevel optimization scheme aiming to minimize the reconstruction error \cite{kunisch2013bilevel,reyes2016total}, e.g., %
\begin{equation}\label{eq:bilevel}
\begin{cases}
\min\limits_{\theta} \sum\limits_{i=1}^N\|\hat{u}_i(\theta)-u_i\|_2^2\\
\text{s.t. } \hat{u}_i(\theta) = \argmin\limits_{u}\Dc_{y_i}(\fwd u) + \Rc_\theta(u).
\end{cases}
\end{equation}
Closely related, in \cite{kobler2020total,kobler2022total}, the authors propose total deep variation regularization. They model the regularizer as $\Rc(u) = \sum_{i,j} \Nc_\theta(u)_{i,j}$ where $\Nc_\theta$ is a U-net style network of a certain structure and the sum is taken over all output pixels of the network. In \cite{kobler2020total,kobler2022total} the solution of the variational regularization problem is formulated as a gradient flow instead of an iterative gradient descent scheme, which allows to phrase the training problem as an optimal control problem. In a discretized setting, however, this leads to a formulation similar to \eqref{eq:bilevel}. The authors of \cite{kobler2020total,kobler2022total} include proofs for existence of solutions of the training problem as well as a stability analysis.

A different approach for training is pursued in \cite{lunz2018adversarial,mukherjee2020learned}. There the authors propose an adversarial training strategy according to
\begin{equation}\label{eq:adversarial_reg}
\min\limits_\theta \E_{u\sim\p(u)}\left[ \Nc_\theta(u) \right] - \E_{y\sim\p(y)}\left[ \Nc_\theta(\fwd^\dagger y) \right] + \E_{\hat{u}\sim\p(\hat{u})}\left[ (\nabla \Nc_\theta(\hat{u})-1)_+^2\right]
\end{equation}
where $\fwd^\dagger$ denotes the (possibly regularized) pseudo-inverse of $\fwd$ and $\p(\hat{u})$ denotes a uniform distribution on straight lines between a sample from $\p(u)$ and a sample $\fwd^\dagger y$ for $y\sim\p(y)$ (cf. \Cref{sec:gans}). The first two terms ensure that the regularizer yields small values on true samples from $\p(u)$ and large values on plain pseudo-inverse reconstructions. The third term, on the other hand, enforces Lipschitz-1 continuity of $\Nc_\theta$ which relates the training problem to the Wasserstein-1 distance \cite{lunz2018adversarial}.

Another different training strategy is pursued in \cite{li2020nett,antholzer2019nett}. Given a data set $(u_i)_i$ of ground truth images, a neural network $\Nc_\theta$ is trained to map $\Nc_\theta(u_i)\approx 0$ for some $i$ and $\Nc_\theta(\fwd^\dagger \fwd u_i)\approx \fwd^\dagger \fwd u_i-u_i$ for the remaining samples. That is, the network maps ground truth images to zero and reconstructions obtained via the pseudo inverse to their residuals to the ground truth. In \cite{li2020nett} the network $\Nc_\theta$ is designed as an auto-encoder $\Nc_\theta = D_{\theta_2}\circ E_{\theta_1}$ and as a regularizer the authors use the norm of the latent variable $\Rc(u) = \|E_{\theta_1}(u)\|_q^q$. In \cite{antholzer2019nett} the regularizer is directly chosen as $\Rc(u) = \|\Nc_\theta(u)\|_2^2$. %

\section{Untrained Models}\label{sec:untrained}
In contrast to supervised and unsupervised methods, where the training of models using training data is an essential building block, untrained neural-network based approaches do not require training, but rather rely on the architecture of the neural network itself to act as a prior for image data.
In this section, we summarize some of the main techniques used in this context. For further information we refer to the review paper \cite{qayyum2022untrained}, which focuses exclusively on untrained neural network priors.

The research field of untrained neural network approaches for inverse problems in imaging was initiated by \cite{ulyanov2018deep}. The main idea of such approaches is to constrain the unknown image data to be in the range of a neural network, whose parameters are optimized to fit the measurement data. This results in approaches of the form
\begin{equation} \label{eq:deep_image_prior_general_objective}
\min_{u,\theta,z} \Dc_y(\fwd u) \qquad \text{s.t. } u=\Nc_\theta(z), 
\end{equation}
which can equivalently be written as
\[
\min_{\theta,z} \Dc_y(\fwd \Nc_\theta(z)).
\]
While the baseline approach above optimizes over both the latent variable $z$ and the network parameters $\theta$, in the original work \cite{ulyanov2018deep} the latent variable $z$ is fixed and initialized with uniform noise between $0$ and $0.1$ as default. The network $\Nc_\theta$ used in \cite{ulyanov2018deep} has a fully convolutional architecture of U-net-type \cite{ronneberger2015u}, with the input and output dimensions being equal. The default network architecture uses skip connections, LeakyReLU as activation function, strides for downsampling and bilinear interpolation for upsampling.

The network used in \cite{ulyanov2018deep} is highly overparametrized (the default architecture has around two million parameters), such that in practice it is capable of fitting also noisy data. However, the authors of \cite{ulyanov2018deep} argue that, following the trajectory of a gradient decent algorithm to minimize the objective \eqref{eq:deep_image_prior_general_objective}, clean images are approximated first and noise or image artifacts appear only in later iterations. Consequently, early stopping is used to achieve a regularization effect.

A second influential work in the context of untrained approaches uses a different strategy to avoid the approximation of noise:  The deep decoder network of \cite{heckel2018deep} uses a simple, non-convolutional architecture with less weight parameters than degrees of freedom in the output. Specifically, for a tensor $z_i \in \R^{n_i \times k_i}$, the tensor $z_{i+1} \in \R^{n_{i+1} \times k_{i+1}}$ of the $(i+1)$ layer for $i=0,\ldots,L-1$ is obtained as
\[ z_{i+1} = \text{cn}\relu(U_i z_i W_i),
\]
where $W_i \in \R^{k_i\times k_{i+1}}$ are the weights, $U_i \in \R^{n_{i+1} \times n_{i}}$ is an upsampling matrix realizing bilinear upsampling (and the identity in case of  $U_{L-1}$) and $\text{cn}$ performs a pointwise normalization of each channel. The final output of the network is then computed as
\[ u = \sigmoid(z_L W_L) .\]
Note that by using the matrix multiplication $z_i W_i$ as weighted connection between consecutive layers, the deep decoder differs significantly from a fully connected neural network in that, for a fixed input and output channel, the same weight is used for all spatial positions. This correspond to what is often referred to as $1\times1$ convolution, and in the continuous setting corresponds to a  convolution with a weighted delta peak at $0$.

Default network dimension used in \cite{heckel2018deep} are $L=6$, $n_0 = 16 \times 16$, $n_d =512 \times 512$ and $k_i = k = 64 $ or $k_i = k = 128$. Like in \cite{ulyanov2018deep}, the input variable $z=z_0$ is chosen randomly and fixed.

The main difference of this model to the original deep image prior of \cite{ulyanov2018deep} is that, due to underparametrization, the deep decoder cannot fit arbitrary data such that no regularization by early stopping is required. This is argued by the authors of \cite{heckel2018deep} by providing in \cite[Proposition 1]{heckel2018deep} a lower bound (that holds with high probability) of how well a single layer deep decoder can approximate zero mean Gaussian noise. 

Again aiming to avoid the approximation of noise, a third class of untrained, neural network based approaches is considered in \cite{habring2022generative}, where also an analysis in the infinite dimensional setting is carried out. There, a variational approach of the form 
\begin{equation} \label{eq:genreg_approach}
\min_{u,v} \Dc_y(\fwd u)  + \lambda_1 \Rc(u-v) + \lambda_2\Gc(v)
\end{equation}
is introduced, where (as special case of the more general framework analyzed in \cite{habring2022generative})

\begin{equation} \label{eq:genreg_prior}
 \Gc(v) = \inf_{\theta,(z_i)_{i=1}^L} \sum_{i=1}^L \|z_i\|_{\Mc} \quad \text{s.t. } \begin{cases}
 v = \Nc_{\theta_1} \circ \ldots \circ \Nc_{\theta_L}(z_L), \\
 z_{i-1} = \Nc_{\theta_i}(z_{i}) \text{ for }i=2,\ldots,L, \\
  \|\theta\|_2 \leq 1, \, M \theta = 0.
 \end{cases}
\end{equation}
Here, $\Nc_{\theta_1} \circ \ldots \circ \Nc_{\theta_L}$ maps the latent variable $z_L$ to image data, and the $\Nc_{\theta_i}$ are convolutional layers mapping latent variables $z_i$ (modeled as Radon measures) via a convolution with filters $\theta_i $ (modeled as $L^2$ functions) to latent variables $z_{i-1}$. Different to the previous two approaches, the layers $\Nc_{\theta_i}$ do not include nonlinearities directly, but a non-linearity is included indirectly via penalizing the Radon norm $\|\cdot \|_{\Mc}$ of the latent variables $z_i$, also at intermediate layers. In addition, for regularization, the weights $\theta_1,\ldots,\theta_L$ are constrained in their $L^2$ norm, and an optional linear constraint (such as having zero mean) is included in the model. 

Another particularity of this approach is that the neural network-type prior $\Gc$ is employed via the infimal-convolution with a second regularization functional $\Rc$. This yields an energy-optimal, additive decomposition of the unknown as $u=(u-v) + v$, where $(u-v)$ is regularized via $\Rc$ and $v$ is regularized via $\Gc$. The reason for this extension is that, as argued in \cite{habring2022generative}, neural-network priors are useful to generate highly structured, texture-like images, while for piecewise smooth images, classical priors such as the total variation \cite{rudin1992tv_mh} or higher-order extensions \cite{bredies2010tgv,holler20ip_review_mh} are already very well suited. This argument is supported in \cite{habring2022generative,habring2023note} also by the analytic result that generator-networks such as \eqref{eq:genreg_prior} can, in the infinite resolution limit, only generate continuous functions as images, making it necessary to combine such priors with other approaches via infimal convolution in order to obtain a regularization that matches the established model of images having sharp edges in the form of jump discontinuities. 

\begin{figure}
\includegraphics[height = 3cm, trim = 11cm 2cm 1cm 5cm, clip]{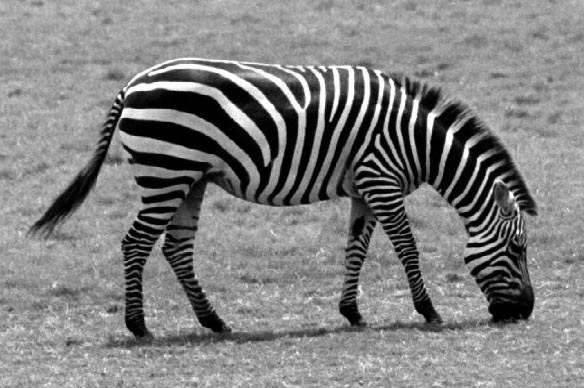}%
\includegraphics[height = 3cm, trim = 8cm 1.4cm 0.5cm 3.6cm, clip]{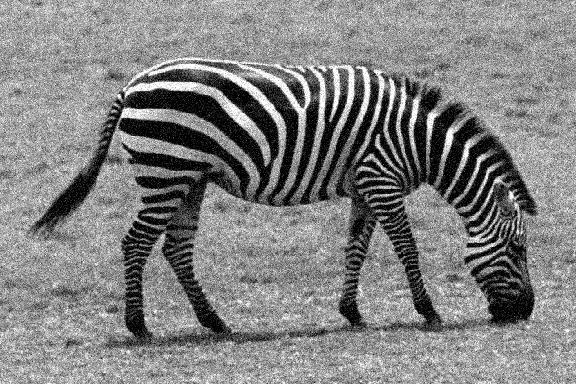}%
\includegraphics[height = 3cm, trim = 11cm 2cm 1cm 5cm, clip]{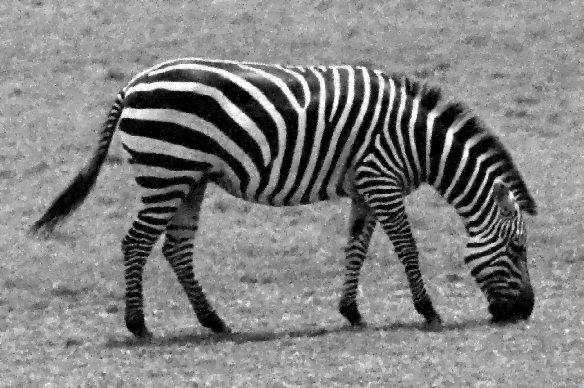}\\
\includegraphics[height = 3.01cm, trim = 11cm 2cm 1cm 5cm, clip]{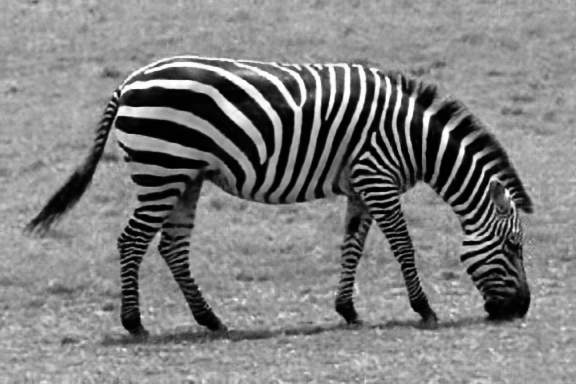}%
\includegraphics[height = 3.01cm, trim = 8cm 1.4cm 0.5cm 3.6cm, clip]{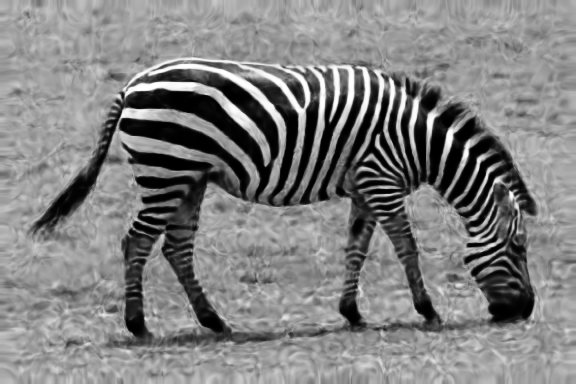}%
\includegraphics[height = 3.01cm, trim = 11cm 2cm 1cm 5cm, clip]{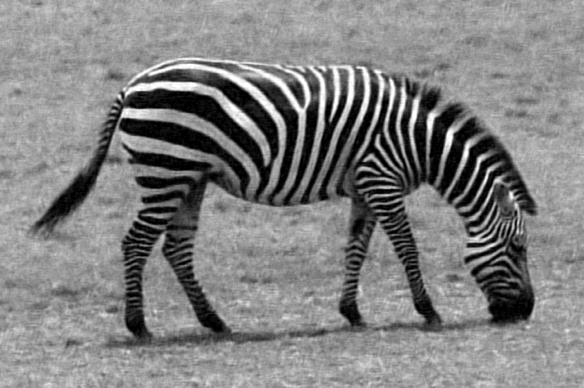}%
\caption{\label{fig:dip_example}Results for denoising with additive Gaussian noise with zero mean and standard deviation 0.1 times the image range. First row: Original, corrupted, and TGV reconstruction \cite{bredies2010tgv}. Second row: Deep image prior \cite{ulyanov2018deep}, deep decoder \cite{heckel2018deep}, and the method from \cite{habring2022generative}.}
\end{figure}

Besides numerical experiments, \cite{habring2022generative} also provides a function-space analysis of the approach \eqref{eq:genreg_approach}, showing in particular well-posedness, a consistency result and the above-mentioned result on regularity of solutions. For an exemplary comparison of the deep image prior, the deep decoder and the approach \eqref{eq:genreg_approach} for denoising we refer to \Cref{fig:dip_example}.

Before moving to a discussion of works on the analysis and extension of deep-image-prior-based approaches for imaging, we also mention that such techniques have also been applied to dynamic data. See for instance \cite{yoo2021time,darestani2021accelerated} for the application of the deep image prior to dynamic MRI or \cite{abdullah2023latent} for the application of an untrained neural-network based approach to the isolation of different types of motion in MRI.
\renewcommand{\ll}{11cm}
\newcommand{\bb}{2}
\newcommand{\rr}{0.5}
\renewcommand{\tt}{5}

\subsection{Regularization Properties of Untrained Neural Networks}\label{sec:reg_untrained}
The, at first glance surprising, capability of the deep image prior to have a regularizing effect in image reconstruction has triggered many works that further analyze and aim to explain this capability. 

The work \cite{heckel2019denoising} analyzes the denoising capabilities of gradient descent with early stopping to fit the deep decoder of \cite{heckel2018deep} to noisy signals. Analytically this is done for a single-layer variant and under the assumption that the ground truth signal is in the span of $p$ trigonometric basis functions, in which case an upper bound (that holds with high probability) on the difference between the denoised and the ground-truth signal is obtained that is of the order $\sigma \sqrt{p/n}$, with $\sigma^ 2 $ the variance of the noise and $n$ the dimension of the signal. In a follow up work \cite{heckel2020compressive}, it is further shown for compressed sensing problems with a Gaussian measurement matrix $\fwd$ that the limit of gradient descent for 
\[ \min_w \frac{1}{2} \|\fwd \Nc (w) - \fwd u^*\|_2 ^2 ,\]
where $\Nc(w)$ is again the output of a single-layer variant of the deep decoder, approximates $u^*$ with high probability.

A conclusion drawn in \cite{heckel2019denoising} is that denoising capabilities of untrained models are primarily attributed to convolutions that use fixed, non-trainable interpolation kernels. Those are usually included in networks via upsampling operators. A similar conclusion is also made in \cite{Liu_2023_ICCV}, which considers spectral effects of upsampling and provides experimental results to conclude that upsampling is key to the denoising effects of untrained models (and, from that, also derives strategies for identifying suitable architectures for the deep image prior).

A different approach that allows to draw conclusions in particular on the regularization properties of untrained models is pursued in \cite{habring2023note}, which was already mentioned above. Here, the authors analyze the regularity of images generated by neural networks. Specifically, they introduce neural networks as mappings between function spaces and show that image-representing functions generated by CNNs with at least two layers are always continuous, and even continuously differentiable if the network input (i.e., the latent variable of the deepest layer) is of bounded variation. A connection to discretized CNNs mapping to images on pixel grids is made via $\Gamma$-convergence, showing that the discretized networks and energies used for training converge to the infinite resolution limit in an appropriate sense. A specific ingredient for the regularity result of \cite{habring2023note} is that the convolution kernels are $L^2$ functions. In the context of network training or fitting an untrained network to data, this is implied by $L ^2$ regularization of the network weights. Hence, a conclusion that can be drawn from \cite{habring2023note} is that $L ^2$ regularization of network weights in the deep image prior approach implies regularity of the resulting, network-generated image and thus avoids the approximation of noise, but on the other hand, hinders the reconstruction of sharp edges in the form of jump discontinuities.

Different to the above approaches, the works \cite{dittmer2020regularization,arndt2022regularization} consider methods related to the deep image prior from a regularization-theoretic perspective. Specifically, they consider the minimization problem
\begin{equation} \label{eq:reg_architecture_main_problem}
\min _B \frac{1}{2}\|\fwd u(B) - y \|_2^2 \quad \text{s.t. } u(B) = \argmin_u \frac{1}{2}\|Bu-y\|_2^2 + \lambda \Rc(u), 
\end{equation} 
with $B$ a linear operator and $\Rc$ a regularization functional, and the solution $u(B)$ as surrogate for $\Nc_\theta(z)$ in the deep image prior objective  \eqref{eq:deep_image_prior_general_objective}. This different formulation is motivated by the fact that, in case $\Rc$ admits a proximal mapping, the iteration
\[ u^{k+1} = \prox _{\tau \lambda \Rc} (u^k - \tau B^*(Bu^k - y)) ,\]
which converges to $u(B)$, resembles, for $L$ steps, the forward pass of a neural network with particular architecture (cf. \Cref{sec:unrolling}). The authors of \cite{dittmer2020regularization} show that, in case $\Rc = \frac{1}{2}\|\cdot \|_2 ^2$ and $B$ is constrained to have the same singular vectors as $\fwd$, the mapping $y \mapsto u(B)$, indeed, corresponds to a regularization method that carries out a filtering of the singular values of the forward map to approximate its pseudo inverse in a way lying in between well established truncated SVD and Tikhonov regularization. The work \cite{arndt2022regularization} further analyzes problems of the form \eqref{eq:reg_architecture_main_problem} and shows in particular that, again if  $\Rc = \frac{1}{2}\|\cdot \|_2 ^2$, they are in fact equivalent to classical Ivanov regularization \cite{Ivanov65_ivanov_regularization}.

\subsection{Exact Recovery of Untrained Network Approaches in Compressed Sensing}
A question very much related to regularization properties of untrained neural-network based approaches is the one of exact recovery capabilities of such approaches in a compressed sensing context.

An early work here is \cite{van2018compressed}, which essentially shows that gradient descent with respect to the network weights $\theta$ on 
\[ \min_\theta \|y-\fwd \Nc_\theta(z) \|_2 ^2 \]
can fit any signal $y$ arbitrarily close, and regards this as evidence that indeed early stopping or additional regularization is required for untrained models. More specifically, this result of  \cite{van2018compressed} is shown for a single-layer network of the form $z \mapsto	\Nc_\theta(z)   = V \cdot \relu(Wz)$ with $\theta=W$ being the unknown parameters and $V$ being a fixed random matrix, and for an orthonormal measurement matrix $\fwd$ and the width of the hidden layer being sufficiently large. In order to avoid overfitting of noisy data, \cite{van2018compressed} further proposes to use additional $\TV$ regularization of the network output and Gaussian-prior-based regularization on the network weights, and also introduces a scheme to learn the parameters of these regularization functionals from data.

In  \cite{heckel2019regularizing}, building on results of \cite{bora2017compressed}, it is shown  that, if $\fwd$ is a random Gaussian matrix with sufficiently many measurements, $y=\fwd u + \zeta$, and if $\hat{\theta}$ minimizes $\theta \mapsto \|y-\fwd\Nc_{\theta} (z)\|_2$ up to an additive $\epsilon$ of the optimum with $z$ fixed, then, with high probability,
\[ \|u - \Nc_ {\hat \theta}(z) \|_2 \leq 6 \min_\theta \|u - \Nc_\theta(z) \|_2 + 3 \|\zeta \|_2 + 2 \epsilon	 + 2 \delta,\]
where $\delta$ is a parameter inversely proportional to the number of measurements.

Another work in this context is \cite{jagatap2019algorithmic}, which proves exact recovery results for linear inverse problems and compressive phase retrieval. Specifically, \cite{jagatap2019algorithmic} introduces algorithms for these two problems that essentially alternate between a gradient step towards data fidelity and a projection to the range of a deep network, and show that the iterates converge to the ground truth in both cases. Main assumptions to achieve this are that the ground truth lies in the range of a generator network, that the sampling matrix obeys a (generalized) restricted isometry property and, in the case of phase retrieval, that the network weights are initialized appropriately and the number of measurements is sufficiently large.

In a more recent work \cite{buskulic:hal-04059168} the gradient flow
\begin{equation} 
\left \{
\begin{aligned}
\dot{\theta}(t) &= - \nabla _\theta \Dc_y ( \fwd (\Nc_\theta(z))) \\
\theta (0) & = \theta_0
\end{aligned}
\right.
\end{equation}
is considered as surrogate for first-order decent algorithms applied to the deep image prior. Besides well-posedness, main results of \cite{buskulic:hal-04059168} are that the state $\theta(t)$ converges to $\theta_\infty$ being a global minimizerer of $\theta \mapsto \Dc_y (\fwd ( \Nc_\theta(z)))$ and, under a restricted injectivity condition, the explicit estimate 
\[ \|\Nc_{\theta(t)}(z) - u\| _2\leq a(t) + b + c ,\] where  $u$ is the true solution such that $y = \fwd(u) + \zeta$ and $\fwd$ is the (potentially non-linear) forward model. Here, $a(t)$ reflects the distance of the current state $\Nc_{\theta(t)}(z)$ to its limit $u_\infty$ (and converges to zero for $t\rightarrow \infty$), $b$ quantifies how close the range of the generator is to the true solution $u$, and $c$ quantifies the level of the noise $\epsilon$. In the linear case, the main assumption for the later result is a restricted injectivity condition on the forward model, which corresponds to the so called minimum conic singular value of $\fwd$ being positive. Using this,  a connection to compressed sensing results like \cite{jagatap2019algorithmic} is drawn in \cite{buskulic:hal-04059168} by showing that, for certain sub-Gaussian measurement matrices, the minimum conic singular value is positive with high probability depending on the number of measurements.

\subsection{Combining Untrained Approaches with other Regularization Functionals}

Based on the works discussed in the previous two subsections, there is a general agreement that overparametrized, untrained neural network approaches such as the deep image prior require some kind of regularization technique to avoid overfitting. While some of the previously discussed works already draw conclusions about suitable regularization strategies as consequence of their results, we now discuss works that explicitly focus on the combination of untrained models with additional regularization functionals.

The work \cite{cheng2019bayesian}, for instance, considers a Bayesian perspective on the deep image prior and uses $\ell^2$ regularization as prior for the network weights and the input. Exploiting the Bayesian perspective, \cite{cheng2019bayesian} further computes the MMSE (instead of the MAP) of posterior
\[
(\theta,z) \mapsto \|y - \Nc_\theta(z)\|_2^2  + \|(\theta,z)\|_2 ^2
\]
for inference. Results of \cite{cheng2019bayesian} indicate that this avoids overfitting and thus the need for early stopping.

A different strategy, pursued for instance in \cite{liu2019image}, is to add $\TV$ regularization on the network output to avoid the approximation of noise, that is, to solve
\[ \min _\theta \|y - \Nc_\theta(z) \|_2 ^2 + \lambda \TV(\Nc_\theta(z))
\]
instead of \eqref{eq:deep_image_prior_general_objective}. Instead of a classical regularization functional like TV, an alternative approach is also to combine the DIP with a trained regularization functional. A popular example in this direction is the combination with regularization by denoising (RED) (cf. \Cref{sec:red}) as proposed in \cite{mataev2019deepred}. Here, the minimization problem is given as
\[ \min _{(\theta,u,z)} \|y - \fwd\Nc_\theta(z) \|_2 ^2 + \lambda u^T (u- D(u)) \qquad \text{s.t. } u = \Nc_\theta(z),
\]
where $D$ can be any pre-trained denoiser. An advantage of this method (and the orginal RED approach) is that, by using specific algorithmic reformulation, approximate solutions to the above optimization problem can be computed without differentiating the learned denoiser.

\section{Function Space Analysis of Neural-network-based Approaches} \label{sec:function_space}

Classical regularization theory for inverse problems usually considers an infinite dimensional function space setting \cite{engl1996regularization}. Beyond the independence of particular discretization schemes, reasons for the later are that crucial concepts of regularization theory for inverse problems, such as ill-posedness of the forward map or the regularity of solutions, are much clearer defined in function spaces. This is in particular also true for inverse problems in imaging, where for instance the space of functions of bounded variation \cite{bvfunctions} provides a rich theory for regularity of image-representing functions, and spaces of measures are the natural setting for grid-independent sparse reconstruction problems.

In the context of neural-network-based approaches for inverse problems in imaging, function space analysis is still a niche topic, which has, however, grown in the past few years. While there exist already many works that analyze neural networks as mappings between function spaces, see for instance \cite{chen1995universal,bach2017breaking,korolev2022two,lanthaler2022error,habring2023note} or \cite{hagemann2023multilevel} for a recent, comprehensive review on infinite dimensional generative models, papers dealing specifically with the analysis of neural-network based approaches for inverse problems are still comparably rare. In this section, we strive to provide a brief overview of works in this context. In addition we refer to \cite{mukherjee2023learned} for a review paper with an emphasis also on the analysis of neural-network based approaches for inverse problems.

An early work dealing with classical regularization theory for neural-network-based approaches to inverse problems is  \cite{schwab2019deep}. Here, the authors consider a specific type of neural network used for post-processing (see \Cref{sec:post-processing} for this type of approaches), for which they provide theoretical results. Given the forward model $\fwd:\Xc \rightarrow \Yc$ and any image-based neural-network $\Nc_\theta:\Xc \rightarrow \Xc$, \cite{schwab2019deep} defines a deep null space network as $L = \id_{\Xc} + P_{\ker(A)} \Nc_\theta:\Xc \rightarrow \Xc$ with $\id_{\Xc}$ the identity operator on $\Xc$. This restricts the image-based regularization of the neural network to the kernel of $\fwd$, such that $\fwd Lu = \fwd u$. Combining this approach with any family $(B_\alpha)_\alpha$ of functionals comprising a regularization method, and assuming Lipschitz continuity of the network $\Nc_\theta:\Xc \rightarrow \Xc$, \cite{schwab2019deep} shows that $(L\circ B_\alpha)$ also comprises a convergent regularization method. For  further extensions and analysis of this approach we refer to \cite{schwab2020big}.

Another, rather influential work providing function space results is \cite{li2020nett}. Here the authors consider a regularizer of the form
\[ \Rc(u) = \psi(\Nc_\theta(u)), \]
called NETT regularizer, where $\psi:\Xc \rightarrow [0,\infty]$ is a given, scalar valued functional and 
\[
\Nc_\theta(u) = (\sigma_L \circ \Theta_L \circ \ldots \sigma_1 \circ \Theta_1) (u),
\]
with $\Theta _i: \Xc_{i-1} \rightarrow \Xc_{i-1/2}$ affine mappings (assumed to be parametrized by $\theta$) and $\sigma_i:\Xc_{i-1/2} \rightarrow \Xc_i$ functions modeling the activation functions. Under assumptions such as weak continuity of the involved mappings and coercivity of $\Rc$, \cite{li2020nett} obtains classical results of regularization theory, in particular well-posedness, convergence and convergence rates under source conditions. 

The follow-up-work \cite{obmann2021augmented} aims to avoid the coercivity assumption of \cite{li2020nett} by considering regularization functionals of the form
\[ \Rc(u) = \psi(\Nc_\theta(u))  + \frac{c}{2}\| u - D_\theta \circ \Nc_\theta(u)\|_2 ^2,\] 
where $\Nc_\theta : \Xc \rightarrow \ell^2$ and $D_\theta : \ell^2 \rightarrow \Xc$ are assumed to be encoder and decoder networks, respectively. Assuming weak sequential continuity in particular of the involved encoder and decoder network, \cite{obmann2021augmented} obtains again classical well-posedness and convergence results, where coercivity is now a consequence of the particular form of $\Rc$. The recent work \cite{haltmeier2023data} further extends this regularization approach to Morozov regularization.

A different, generator (or synthesis) based approach towards neural-network-based regularization in function space is considered in \cite{obmann2020deep}. There, the authors consider a variational regularization approach of the form
\[
\min_z \| \fwd ( \Nc_{\theta,\lambda} (z)) - y \|_2^2 + \lambda \|z\|_{1,w},
\]
where $\Nc_{\theta,\lambda}:Z \rightarrow \Xc$ is a generator network (depending also on the regularization parameter $\lambda$) and $\|\cdot \|_{1,w}$ is a weighted $\ell^1$ norm. Assuming the generator network $\Nc_{\theta,\lambda}$ to be weakly sequentially continuous, \cite{obmann2020deep} obtains again well-posedness and convergence results for this approach.

Function space analysis for a related, generator-based regularization strategy was also considered in the aforementioned  \cite{habring2022generative} via a variational approach of the form
\begin{equation} \label{eq:genreg_approach_function_space_section}
\min_{u,v} \Dc_y(\fwd u)  + \lambda_1 \Rc(u-v) + \lambda_2\Gc(v).
\end{equation}
Recall that here, $\Gc$ is a prior based on a generator network as defined in \eqref{eq:genreg_prior}, and $\Rc$ is an optional, classical regularization functional. Different to the approach of \cite{obmann2020deep} above, \cite{habring2022generative} considers an untrained model in the spirit of the deep image prior \cite{ulyanov2018deep}. Modeling latent spaces as Radon measures and filter kernls as $L^2$ functions, \cite{habring2022generative} proves weak-to-weak continuity of the involved, bilinear convolution operators and coercivity of the overall approach, from which well-posedness and convergence results follow using standard techniques. A particularity of \cite{habring2022generative} is that the explicit definition of convolutional neural networks as mapping between function spaces also allows to analyze regularity of solutions of \eqref{eq:genreg_approach_function_space_section}. In particular, in \cite{habring2022generative} it is shown that the output of the generator network is always a continuous function, implying that solutions $u,v$ of \eqref{eq:genreg_approach_function_space_section} satisfy $u = (u-v) + v \in \BV(\Omega) + C(\Omega)$ in case $\Rc $ is the total variation functional and $\Omega \subset \R^d$ is the domain of the image representing functions.

In this context, we also recall the work \cite{habring2023note}, which is already discussed in  \Cref{sec:reg_untrained} and analyses the regularity of images generated by convolutional neural networks. In the context of regularization theory for neural network approaches in function space, it is also worth noting that, as supplementary result, \cite[Theorem 3.15]{habring2023note} provides well-posedness of the training of a large class of neural networks in function space by empirical risk minimization, without posing assumptions such as weak continuity or coercivity on the network.

An approach similar to the NETT regularizer of  \cite{li2020nett} is the adversarial regularization approach of \cite{lunz2018adversarial} (see also \Cref{sec:learned_prior}). Here, the regularization functional is of the form
\[ 
u \mapsto \Nc_\theta(u),
\]
where $\Nc_\theta$ is a network trained to distinguish images coming from the true image disctrubtion and images with artifacts being obtained via a direct reconstruction method. Besides well-posedness results under lower semi-continuity and coercivity assumptions, main results of \cite{lunz2018adversarial} are that gradient descent on a perfectly trained regularizer moves images towards the ground truth image distribution, and that, if the training data lies on a manifold, a perfectly trained model projects to this manifold. 
An extension of \cite{lunz2018adversarial} towards training input-convex adversarial regularizers is provided in \cite{mukherjee2020learned}, where input-convex neural networks in function spaces are defined as integral transforms with partially non-negative kernels to obtain convexity.

The work \cite{alberti2021learning} analyses the learning of an optimal Tikhonov regularization approach in Hilbert spaces from a statistical perspective, where the unknowns $u$ and the data $y = \fwd u + \zeta$ are modeled as random variables in Hilbert spaces with joint distribution $\p(u,y)$. Omitting - very interesting - technicalities from \cite{alberti2021learning} that are necessary to rigorously define the statistical setup in a function space setting, the main idea is to solve
\begin{equation} \label{eq:alberti_optimal_tikhonov}
\min_{h,B} \E_{(u,y)\sim p(u,y)} \|R_{h,B}(y) - u \|_{\Xc}^2 \quad \text{where } R_{h,B}(y) = \argmin_u \Dc_y(\fwd u) + \|B^{-1} (u-h)\|_{\Xc} ^2.
\end{equation}
As shown in \cite{alberti2021learning}, $(h,B) = (\mu,\Sigma_x^{1/2})$, with $\mu$ and $\Sigma_x$ the mean and covariance of the random variable describing the ground truth image distribution, is always a solution to \eqref{eq:alberti_optimal_tikhonov}, which is in particular independent of the forward model. These results are further extended in \cite{alberti2021learning} to training with finitely many samples via convergence results.

Related to the previous work, in \cite{kabri2022convergent} the authors consider the learning of a filter function either for spectral regularization or for the filtered-backprojection for the Radon transform. In the former case, assuming the forward operator $\fwd$ to be compact and to be given via the singular value expansion as
\[ \fwd u = \sum_{n=1}^\infty \sigma_n \langle u,u_n\rangle v_n ,\]
the spectral regularization approach takes the form
\[ R(y^\delta,g_\delta) = \sum_{n=1}^\infty g_\delta(\sigma_n) \langle y^\delta ,v_n\rangle u_n ,\]
where $g_\delta$ is a learned function with $g_\delta(\sigma) \rightarrow	1/\sigma $ as $\delta \rightarrow 0$ and $\delta$ is the noise level. A result of \cite{kabri2022convergent} is that, in the ideal case of minimizing the true risk, the optimal filter function is given as
\[ g(\sigma_n) = \frac{\sigma_n \Pi_n}{\sigma_n^2 \Pi_n + \Delta_n} \]
where $\Pi_n = \E_u (\langle u,u_n\rangle^2 )$ and $\Delta_n = \E_\zeta (\langle \zeta,v_n\rangle^2 )$, with $\zeta$ the additive measurement noise. Beyond that, \cite{kabri2022convergent} provides convergence results for the learned regularization.

The work \cite{arndt2023invertible} considers a residual network, i.e., networks of the form
\[ \varphi_\theta(x) = x - \Nc_\theta(x) \]
to approximate the forward operator $\fwd$, and then exploits inversion of $\varphi_\theta$ via a fixed point iteration to approximate $\fwd^{-1}$. Following this approach, \cite{arndt2023invertible}  provides well-posedness and convergence results, and also considers particular architectures of $\Nc_\theta$ and the relation to spectral regularization methods.

In a different context, \cite{ebner2022plug} considers the analysis of PnP methods (see \Cref{sec:pnp} for details) as regularization methods. Specifically, they consider iterations of the form
\[ u_{n+1} = (\Nc_\theta(\lambda,\cdot) \circ (\id - s\nabla _u \Dc_{y}) ) (u^n) ,\]
where $\Nc_\theta(\lambda,\cdot)$ is a learned denoiser depending on the regularization parameter $\lambda$. Assuming that that $ \Nc_\theta(\lambda,\cdot)$ is a contraction and that $ \Nc_\theta(\lambda,\cdot)$ converges to the identity as $\lambda \rightarrow 0$ in an appropriate sense, \cite{ebner2022plug} shows that the mapping of data to fixed points of the above regularization is a convergent regularization method.

The work \cite{aspri2020data} introduces data driven regularization via including a projection to training data. That this, given training pairs $u_i$ and $y_i = \fwd u_i$ for $i=1,\ldots,n$ as training data, \cite{aspri2020data} analyzes the regularization properties of the generalized inverse of $\fwd P_{U_n}$, with $P_{U_n} $ a projection to $U_n = \vspan((u_i)_{i=1} ^n)$, which is shown to be equal to 
\[ \fwd^{-1} P_{Y_n} .\]
In \cite{aspri2020data} the explicit expression of $A^{-1} P_{Y_n}$ is further investigated in terms of orthonormal bases obtained from the training data, stability of the orthonormalization procedure and regualrization properties of the overall resulting method. In particular, convergence in the case of noise free and noisy data is shown under appropriate conditions.

The paper \cite{aspri2020datad} analyzes a data-driven Landweber iteration for non-linear inverse problems. The iteration is given as
\[ u^{k+1} = u^k - \fwd^\prime (u^k) ^*( \fwd(u^k) - y) - \lambda_k B^\prime (u^k)^* (B(u^k) - y ) ,
\]
where $B:\Xc \rightarrow \Yc$ is a learned operator satisfying $B(u_i) = y_i$ for some training data $(u_i,y_i)_{i=1}^N$. The motivation for this data-driven adaption is that the damping term appearing in the classical Landweber iteration can be interpred as a prior on the solution, and hence the penalization of $\|B(u) - y\|_2$ allows to include additional expert knowledge in the form of training data. Under appropriate assumptions, in particular the tangential cone condition, \cite{aspri2020datad}  obtains convergence results for the proposed modified Landweber iteration. Further results on convergence and stability of this method are also obtained in \cite{aspri2023analysis}, and in the work \cite{gao2023fast} the method is further extended with a convex penalty term and acceleration techniques.

The work \cite{scherzer2023gauss} analyzes a Gauss-Newton method in a setting conceptually related to sparse coding techniques and untrained models like the deep image prior. That is,  \cite{scherzer2023gauss} analyzes the Gauss-Newton method for solving 
\[ \fwd \Nc_\theta(\;\cdot\;) = y \]
with respect to $\theta$, where $\theta \mapsto \Nc_\theta$ maps finite dimensional parameters to a shallow neural network representing the unknown as function in $L^2$. A result of \cite{scherzer2023gauss}  is in particular local quadratic convergence of the resulting method under appropriate conditions, in particular linear independence of activation functions and their derivatives.

Related to this is also the work \cite{alberti2022continuous}, which considers pretrained deep generative models $\Nc_\theta : \R^S \rightarrow \Xc$ mapping from a finite dimensional latent space to image space. The work \cite{alberti2022continuous} rigorously introduces deep generative models mapping to function space via a wavelet-based multi-resolution analysis. It further investigates injectivity of the generator and provides a stability result for inverting non-linear forward models on the range space of the generator.

The work \cite{baldassari2023conditional} further extends score based generative models to infinite dimensions (see also \cite{hagemann2023multilevel} for different approaches towards infinite dimensional generative models), and employs them to sample from a conditional distribution given the observed data in a Bayesian linear inverse problem setting.%

\section{Quantifying Uncertainty in Learned Models} \label{sec:uncertainty_quantification}
Uncertainty quantification is a central topic both in machine learning (see for instance \cite{abdar2021review} for a recent review) and inverse problems (see \cite{tenorio2017introduction}). For the latter, recent works commonly adapt the Bayesian perspective on this topic \cite{stuart2010inverse}, and we focus on this viewpoint also in this section. 

The general approach in Bayesian uncertainty quantification is rather simple: As explained in \Cref{sec:introduction}, describing the data at hand as random variables $(u,y)$ with a joint distribution $\p(u,y)$, where typically $y = \fwd(u) + \zeta$ with the random variable $\zeta $ describing measurement noise, Bayes' theorem provides the posterior distribution of $u$ given an observation $y$ as
\[ \p(u|y) = \frac{\p(y|u)\p(u)}{\p(y)} \propto \p(y|u)\p(u) .\]
Thus, in theory, given the data likelihood $\p(y|u)$ and the prior $\p(u)$, a full statistical discription of the quantity of interest is available, from which point estimates (such as the MAP and MMSE) and information on uncertainty (such as variance or error quantiles) can be obtained. In practice, the situation is of course more difficult, and all main aspects of this approach comprise important research topics: Defining the data likelihood $\p(y|u)$ is mostly a question of obtaining appropriate forward and noise models (and is often not in the focus of Bayesian approaches to inverse problems in imaging). Obtaining the prior $\p(u)$ has received a lot of attention of the research community, where research has moved from hand-crafted priors to priors that are explicitly or implicitly learned from data, which was the main content of \Cref{sec:unsupervised}. Another vivid topic of research are methods that allow for efficient sampling from the posterior, with MCMC methods being often in the focus of attention (see for instance \cite{durmus19langevin,durmus2019analysis,durmus2022proximal,laumont2022bayesian,habring2023subgradient}).

Machine learning models that either learn the prior explicitly (see \Cref{sec:learned_prior}) or implicitly (see \Cref{sec:learned_generators}) can usually directly be extended towards uncertainty quantification, at least if sampling is computationally feasible. Uncertainty quantification in this context typically means to compute pixel-wise marginal variances and regard them as indicator for uncertainty, see for instance \cite{zach2023stable,zach2022computed}.
Generator-based priors, i.e., priors where the image manifold is described via a network $\Nc_\theta(z)$ with trained parameters and the latent variable $z$ following a specified prior distribution, are analyzed in a Bayesian context for example in \cite{holden2022bayesian} (see also \cite{gonzales2022solving} for a related approach which, however, focuses on robust MAP estimation). The work \cite{holden2022bayesian} shows well-posedness of the posterior resulting from the incorporation of a learned, generator-based prior and also the existence of posterior moments, which can be sampled with MCMC techniques. Exploiting low dimensionality of the latent space, uncertainty can be visualized in this context via considering the covariance matrix for the latent variable, and mapping forward information on the principal components of this matrix to image space.

In \cite{narnhofer2022bayesian} the authors impose a prior distribution $\p(\theta)$ on the parameters $\theta$ of a learned regularizer $\Rc(u) = \Nc_\theta(u)$. That is, instead of estimating $\theta$ directly, in a Bayesian fashion they try to infer the distribution of $\theta$. Subsequently, this allows to quantify the uncertainty with respect to the optimal parameters $\theta$. The authors provide experiments with the TDV prior \cite{kobler2020total}.

The work \cite{schlemper2018bayesian} uses Markov chain dropout to derive uncertainty estimates for deep-learning-based MRI reconstruction, and also qualitatively relates uncertainty to the reconstruction error.
The work \cite{ramzi2020denoising} uses a denoising autoencoder to approximate the gradient of the prior distribution in a Bayesian setting for MRI imaging, and combines this with sampling techniques to assess the uncertainty of particular features of the reconstructed images.
Also in the context of Bayesian MR imaging, the work \cite{luo2023bayesian} uses diffusion models both for reconstruction and for the computation of uncertainty maps via sampling.

Also plug and play methods were recently extended to allow for uncertainty quantification. In \cite{laumont2022bayesian} for instance, the authors use Tweedies formula to re-interpret the gradient of the log-prior, appearing in sampling algorithms, as MMSE denoiser, which opens the door to again perform sampling of a posterior distribution implicitly described by pre-trained denoisers. The work \cite{laumont2022bayesian}  analyzes this approach, in particular w.r.t. convergence of the resulting algorithms, and uses sampling to compute pixel-wise marginal posterior standard deviation.

A limitation of classical Bayesian uncertainty quantification methods, in particular in connection with learned, deep-neural-network-based priors, is that quantitative guarantees, such as error bounds that hold with high probability, can hardly be rigorously ensured. There are manifold reasons for this: It is hard to know or quantify the discrepancy between the true prior and the one learned on a finite sample; even in the idealized case of sufficient training samples, non-convexity of learning problems with deep neural networks usually prevent guarantees that the global minimum of the loss was found; and even if the prior was trained perfectly, sampling algorithms such as Langevin sampling usually require log-concave densities \cite{durmus2019analysis}, which is not the case for deep-network-based models. 
The works \cite{angelopoulos2022image,narnhofer2022posterior} consider an alternative strategy towards uncertainty quantification for inverse imaging problems, which is not directly affected by these limitations. Building on techniques of conformal prediction from \cite{romano2019conformalized}, the main idea of \cite{angelopoulos2022image,narnhofer2022posterior} is to develop machine-learning-based predictors for error intervals, but to regard them as black-box predictors and use calibration data together with quantile estimates to obtain corrected versions of those predictors, whose error bounds are valid with high probability independent on assumptions on the underlying data distribution. Specifically, \cite{narnhofer2022posterior} considers the following setting: Assume we are given a learned prior and a specific inverse problem. Given an iid \emph{calibration} data set $(u_i,y_i)_i$, we can use the prior to compute point estimates $\hat{u}_i$ of $u_i$ as well as estimates of the posterior variances $v_i=\var[u|y=y_i]$. The labeled data set then allows to estimate the relation between $v_i$ and the error $\hat{u}_i-u_i$. Given a new unseen observation $y_{\text{new}}$ we can compute $v_\text{new}=\var[u|y=y_{\text{new}}]$ and estimate the reconstruction error based on the learned relation between posterior variance and error. In \cite{narnhofer2022posterior} (relying on results of \cite{romano2019conformalized}) it is proven that the obtained error estimates hold with high probability (the probability taken over calibration data set, new observation $y$ and the estimate of $u$), independent of assumptions on the data distribution or sampling algorithms.

\bibliographystyle{plain}
\bibliography{lit_dat}

\end{document}

%% file: notation.tex
\DeclareMathOperator{\Rc}{\mathcal{R}}
\DeclareMathOperator{\Dc}{\mathcal{D}}
\DeclareMathOperator{\Nc}{\mathcal{N}}
\DeclareMathOperator{\Yc}{\mathcal{Y}}
\DeclareMathOperator{\Xc}{\mathcal{X}}
\DeclareMathOperator{\Zc}{\mathcal{Z}}
\DeclareMathOperator{\var}{\mathbb{V}}
\DeclareMathOperator{\Gc}{\mathcal{G}}
\DeclareMathOperator{\KL}{KL}
\renewcommand{\Mc}{\mathcal{M}}
\DeclareMathOperator{\fwd}{A}
\let\p\relax

\DeclareMathOperator{\p}{p}
\DeclareMathOperator{\q}{q}

\DeclareMathOperator{\rg}{Rg}

\DeclareMathOperator{\E}{\mathbb{E}}
\DeclareMathOperator{\R}{\mathbb{R}}
\DeclareMathOperator{\relu}{ReLU}
\DeclareMathOperator{\J}{J}
\DeclareMathOperator{\tr}{tr}
\DeclareMathOperator{\sigmoid}{sigmoid}
\DeclareMathOperator*{\argmax}{arg\,max}
\DeclareMathOperator*{\argmin}{arg\,min}
\renewcommand{\d}{\text{d}}
\DeclareMathOperator{\TV}{TV}
\DeclareMathOperator{\BV}{BV}
\DeclareMathOperator{\id}{I} 
\DeclareMathOperator{\I}{I} 
\DeclareMathOperator{\prox}{prox} 

\DeclareMathOperator{\vspan}{span}

\forestset{
  dir switch tree 2/.style={
    for tree={
      if level=1{
        edge path={
          \noexpand\path [\forestoption{edge}] (!u.parent anchor) -- ++(0,-1.5em) -| (.child anchor)\forestoption{edge label};
        },
        for descendants={
          child anchor=west,
          align=left,
          edge path={
            \noexpand\path [\forestoption{edge}] (!u.parent anchor) -- ++(-.5em,0) |- (.child anchor)\forestoption{edge label};
          },
          before computing xy={
            l=.5em
          }
        },
        before computing xy={
          l=3em
        },
        for tree={
          parent anchor=west,
          anchor=mid west,
          grow'=0,
          font=\sffamily,
          if n children=0{}{
            delay={
              prepend={[,phantom, calign with current]}
            }
          },
        },
      }{},
    },
  },
}
\tikzset{
  colour me/.style={
      draw=colour#1,
      fill=colour#1!20!white,
      text=colour#1!25!black,
  }
}

\colorlet{colour0}{blue}
\colorlet{colour1}{green}
\colorlet{colour2}{red}
\colorlet{colour3}{orange}
\colorlet{colour4}{yellow}
\colorlet{colour5}{magenta}

%% file: paper.bbl
\begin{thebibliography}{100}

\bibitem{abdar2021review}
Moloud Abdar, Farhad Pourpanah, Sadiq Hussain, Dana Rezazadegan, Li~Liu,
  Mohammad Ghavamzadeh, Paul Fieguth, Xiaochun Cao, Abbas Khosravi, U~Rajendra
  Acharya, et~al.
\newblock A review of uncertainty quantification in deep learning: Techniques,
  applications and challenges.
\newblock {\em Information fusion}, 76:243--297, 2021.

\bibitem{abdullah2023latent}
Abdullah Abdullah, Martin Holler, Karl Kunisch, and Malena~Sabate Landman.
\newblock Latent-space disentanglement with untrained generator networks for
  the isolation of different motion types in video data.
\newblock In {\em International Conference on Scale Space and Variational
  Methods in Computer Vision}, pages 326--338. Springer, 2023.

\bibitem{adler2017solving}
Jonas Adler and Ozan Öktem.
\newblock Solving ill-posed inverse problems using iterative deep neural
  networks.
\newblock {\em Inverse Problems}, 33(12):124007, 2017.

\bibitem{adler2018learned}
Jonas Adler and Ozan Öktem.
\newblock Learned primal-dual reconstruction.
\newblock {\em IEEE Transactions on Medical Imaging}, 37(6):1322--1332, 2018.

\bibitem{alain2014regularized}
Guillaume Alain and Yoshua Bengio.
\newblock What regularized auto-encoders learn from the data-generating
  distribution.
\newblock {\em The Journal of Machine Learning Research}, 15(1):3563--3593,
  2014.

\bibitem{alberti2021learning}
Giovanni~S Alberti, Ernesto De~Vito, Matti Lassas, Luca Ratti, and Matteo
  Santacesaria.
\newblock Learning the optimal tikhonov regularizer for inverse problems.
\newblock {\em Advances in Neural Information Processing Systems},
  34:25205--25216, 2021.

\bibitem{alberti2022continuous}
Giovanni~S Alberti, Matteo Santacesaria, and Silvia Sciutto.
\newblock Continuous generative neural networks.
\newblock {\em arXiv preprint arXiv:2205.14627}, 2022.

\bibitem{altekruger2022patchnr}
Fabian Altekr{\"u}ger, Alexander Denker, Paul Hagemann, Johannes Hertrich,
  Peter Maass, and Gabriele Steidl.
\newblock Patchnr: Learning from small data by patch normalizing flow
  regularization.
\newblock {\em arXiv preprint arXiv:2205.12021}, 2022.

\bibitem{bvfunctions}
Luigi Ambrosio, Nicola Fusco, and Diego Pallara.
\newblock {\em Functions of Bounded Variation and Free Discontinuity Problems}.
\newblock Oxford Mathematical Monographs, 2000.

\bibitem{anderson1982reverse}
Brian~D.O. Anderson.
\newblock Reverse-time diffusion equation models.
\newblock {\em Stochastic Processes and their Applications}, 12(3):313--326,
  1982.

\bibitem{angelopoulos2022image}
Anastasios~N Angelopoulos, Amit~Pal Kohli, Stephen Bates, Michael Jordan,
  Jitendra Malik, Thayer Alshaabi, Srigokul Upadhyayula, and Yaniv Romano.
\newblock Image-to-image regression with distribution-free uncertainty
  quantification and applications in imaging.
\newblock In {\em International Conference on Machine Learning}, pages
  717--730. PMLR, 2022.

\bibitem{antholzer2019nett}
Stephan Antholzer, Johannes Schwab, Johnnes Bauer-Marschallinger, Peter
  Burgholzer, and Markus Haltmeier.
\newblock {NETT regularization for compressed sensing photoacoustic
  tomography}.
\newblock In Alexander~A. Oraevsky and Lihong~V. Wang, editors, {\em Photons
  Plus Ultrasound: Imaging and Sensing 2019}, volume 10878, page 108783B.
  International Society for Optics and Photonics, SPIE, 2019.

\bibitem{arjomand2017deep}
Siavash Arjomand~Bigdeli, Matthias Zwicker, Paolo Favaro, and Meiguang Jin.
\newblock Deep mean-shift priors for image restoration.
\newblock In I.~Guyon, U.~Von Luxburg, S.~Bengio, H.~Wallach, R.~Fergus,
  S.~Vishwanathan, and R.~Garnett, editors, {\em Advances in Neural Information
  Processing Systems}, volume~30. Curran Associates, Inc., 2017.

\bibitem{arjovsky2017wasserstein}
Martin Arjovsky, Soumith Chintala, and L{\'e}on Bottou.
\newblock {W}asserstein generative adversarial networks.
\newblock In Doina Precup and Yee~Whye Teh, editors, {\em Proceedings of the
  34th International Conference on Machine Learning}, volume~70 of {\em
  Proceedings of Machine Learning Research}, pages 214--223. PMLR, 2017.

\bibitem{arndt2022regularization}
Clemens Arndt.
\newblock Regularization theory of the analytic deep prior approach.
\newblock {\em Inverse Problems}, 38(11):115005, 2022.

\bibitem{arndt2023invertible}
Clemens Arndt, Alexander Denker, S{\"o}ren Dittmer, Nick Heilenk{\"o}tter,
  Meira Iske, Tobias Kluth, Peter Maass, and Judith Nickel.
\newblock Invertible residual networks in the context of regularization theory
  for linear inverse problems.
\newblock {\em arXiv preprint arXiv:2306.01335}, 2023.

\bibitem{arridge2019solving}
Simon Arridge, Peter Maass, Ozan Öktem, and Carola-Bibiane Schönlieb.
\newblock Solving inverse problems using data-driven models.
\newblock {\em Acta Numerica}, 28:1–174, 2019.

\bibitem{asim2020invertible}
Muhammad Asim, Max Daniels, Oscar Leong, Ali Ahmed, and Paul Hand.
\newblock Invertible generative models for inverse problems: mitigating
  representation error and dataset bias.
\newblock In Hal~Daumé III and Aarti Singh, editors, {\em Proceedings of the
  37th International Conference on Machine Learning}, volume 119 of {\em
  Proceedings of Machine Learning Research}, pages 399--409. PMLR, 2020.

\bibitem{asim2018solving}
Muhammad Asim, Fahad Shamshad, and Ali Ahmed.
\newblock Solving bilinear inverse problems using deep generative priors.
\newblock {\em CoRR, abs/1802.04073}, 3(4):8, 2018.

\bibitem{asim2020blind}
Muhammad Asim, Fahad Shamshad, and Ali Ahmed.
\newblock Blind image deconvolution using deep generative priors.
\newblock {\em IEEE Transactions on Computational Imaging}, 6:1493--1506, 2020.

\bibitem{aspri2020datad}
Andrea Aspri, Sebastian Banert, Ozan {\"O}ktem, and Otmar Scherzer.
\newblock A data-driven iteratively regularized landweber iteration.
\newblock {\em Numerical Functional Analysis and Optimization},
  41(10):1190--1227, 2020.

\bibitem{aspri2020data}
Andrea Aspri, Yury Korolev, and Otmar Scherzer.
\newblock Data driven regularization by projection.
\newblock {\em Inverse Problems}, 36(12):125009, 2020.

\bibitem{aspri2023analysis}
Andrea Aspri and Otmar Scherzer.
\newblock Analysis of generalized iteratively regularized landweber iterations
  driven by data.
\newblock {\em arXiv preprint arXiv:2312.03337}, 2023.

\bibitem{bach2017breaking}
Francis Bach.
\newblock Breaking the curse of dimensionality with convex neural networks.
\newblock {\em The Journal of Machine Learning Research}, 18(1):629--681, 2017.

\bibitem{bai2020deep}
Yanna Bai, Wei Chen, Jie Chen, and Weisi Guo.
\newblock Deep learning methods for solving linear inverse problems: Research
  directions and paradigms.
\newblock {\em Signal Processing}, 177:107729, 2020.

\bibitem{baldassari2023conditional}
Lorenzo Baldassari, Ali Siahkoohi, Josselin Garnier, Knut Solna, and Maarten~V
  de~Hoop.
\newblock Conditional score-based diffusion models for bayesian inference in
  infinite dimensions.
\newblock {\em arXiv preprint arXiv:2305.19147}, 2023.

\bibitem{beck2009fista}
Amir Beck and Marc Teboulle.
\newblock A fast iterative shrinkage-thresholding algorithm for linear inverse
  problems.
\newblock {\em SIAM Journal on Imaging Sciences}, 2(1):183--202, 2009.

\bibitem{berner2021modern}
Julius Berner, Philipp Grohs, Gitta Kutyniok, and Philipp Petersen.
\newblock The modern mathematics of deep learning.
\newblock {\em arXiv preprint arXiv:2105.04026}, pages 86--114, 2021.

\bibitem{bigdeli2017image}
Siavash~Arjomand Bigdeli and Matthias Zwicker.
\newblock Image restoration using autoencoding priors.
\newblock {\em arXiv preprint arXiv:1703.09964}, 2017.

\bibitem{bora2017compressed}
Ashish Bora, Ajil Jalal, Eric Price, and Alexandros~G Dimakis.
\newblock Compressed sensing using generative models.
\newblock In {\em International conference on machine learning}, pages
  537--546. PMLR, 2017.

\bibitem{boyd2011admm}
Stephen Boyd, Neal Parikh, Eric Chu, Borja Peleato, and Jonathan Eckstein.
\newblock Distributed optimization and statistical learning via the alternating
  direction method of multipliers.
\newblock {\em Foundations and Trends in Machine Learning}, 3(1):1--122, 2011.

\bibitem{holler20ip_review_mh}
Kristian Bredies and Martin Holler.
\newblock Higher-order total variation approaches and generalisations.
\newblock {\em Inverse Problems. Topical Review}, 36(12):123001, 2020.

\bibitem{bredies2010tgv}
Kristian Bredies, Karl Kunisch, and Thomas Pock.
\newblock Total generalized variation.
\newblock {\em SIAM Journal on Imaging Sciences}, 3(3):492--526, 2010.

\bibitem{bredies2008linear}
Kristian Bredies and Dirk~A. Lorenz.
\newblock Linear convergence of iterative soft-thresholding.
\newblock {\em Journal of Fourier Analysis and Applications}, 14(5):813--837,
  2008.

\bibitem{brifman2016turning}
Alon Brifman, Yaniv Romano, and Michael Elad.
\newblock Turning a denoiser into a super-resolver using plug and play priors.
\newblock In {\em 2016 IEEE International Conference on Image Processing
  (ICIP)}, pages 1404--1408, 2016.

\bibitem{buades2005nonlocal}
A.~Buades, B.~Coll, and J.-M. Morel.
\newblock A non-local algorithm for image denoising.
\newblock In {\em 2005 IEEE Computer Society Conference on Computer Vision and
  Pattern Recognition (CVPR'05)}, volume~2, pages 60--65 vol. 2, 2005.

\bibitem{buskulic:hal-04059168}
Nathan Buskulic, Jalal Fadili, and Yvain Qu{\'e}au.
\newblock Convergence and recovery guarantees of unsupervised neural networks
  for inverse problems.
\newblock {\em arXiv preprint arXiv:2309.12128}, 2023.

\bibitem{chambolle2011pd}
Antonin Chambolle and Thomas Pock.
\newblock A first-order primal-dual algorithm for convex problems with
  applications to imaging.
\newblock {\em Journal of Mathematical Imaging and Vision}, 40(1):120--145,
  2011.

\bibitem{chan2017plug}
Stanley~H. Chan, Xiran Wang, and Omar~A. Elgendy.
\newblock Plug-and-play admm for image restoration: Fixed-point convergence and
  applications.
\newblock {\em IEEE Transactions on Computational Imaging}, 3(1):84--98, 2017.

\bibitem{chen1997convergence}
George H-G. Chen and R.~T. Rockafellar.
\newblock Convergence rates in forward--backward splitting.
\newblock {\em SIAM Journal on Optimization}, 7(2):421--444, 1997.

\bibitem{chen2017low}
Hu~Chen, Yi~Zhang, Mannudeep~K. Kalra, Feng Lin, Yang Chen, Peixi Liao, Jiliu
  Zhou, and Ge~Wang.
\newblock Low-dose ct with a residual encoder-decoder convolutional neural
  network.
\newblock {\em IEEE Transactions on Medical Imaging}, 36(12):2524--2535, 2017.

\bibitem{chen1995universal}
Tianping Chen and Hong Chen.
\newblock Universal approximation to nonlinear operators by neural networks
  with arbitrary activation functions and its application to dynamical systems.
\newblock {\em IEEE transactions on neural networks}, 6(4):911--917, 1995.

\bibitem{chen2017trainable}
Yunjin Chen and Thomas Pock.
\newblock Trainable nonlinear reaction diffusion: A flexible framework for fast
  and effective image restoration.
\newblock {\em IEEE Transactions on Pattern Analysis and Machine Intelligence},
  39(6):1256--1272, 2017.

\bibitem{cheng2019bayesian}
Zezhou Cheng, Matheus Gadelha, Subhransu Maji, and Daniel Sheldon.
\newblock A bayesian perspective on the deep image prior.
\newblock In {\em Proceedings of the IEEE/CVF Conference on Computer Vision and
  Pattern Recognition}, pages 5443--5451, 2019.

\bibitem{chung2023parallel}
Hyungjin Chung, Jeongsol Kim, Sehui Kim, and Jong~Chul Ye.
\newblock Parallel diffusion models of operator and image for blind inverse
  problems.
\newblock In {\em Proceedings of the IEEE/CVF Conference on Computer Vision and
  Pattern Recognition (CVPR)}, pages 6059--6069, 2023.

\bibitem{chung2022diffusion}
Hyungjin Chung, Jeongsol Kim, Michael~T Mccann, Marc~L Klasky, and Jong~Chul
  Ye.
\newblock Diffusion posterior sampling for general noisy inverse problems.
\newblock {\em arXiv preprint arXiv:2209.14687}, 2022.

\bibitem{chung2022improving}
Hyungjin Chung, Byeongsu Sim, Dohoon Ryu, and Jong~Chul Ye.
\newblock Improving diffusion models for inverse problems using manifold
  constraints.
\newblock In S.~Koyejo, S.~Mohamed, A.~Agarwal, D.~Belgrave, K.~Cho, and A.~Oh,
  editors, {\em Advances in Neural Information Processing Systems}, volume~35,
  pages 25683--25696. Curran Associates, Inc., 2022.

\bibitem{chung2022come}
Hyungjin Chung, Byeongsu Sim, and Jong~Chul Ye.
\newblock Come-closer-diffuse-faster: Accelerating conditional diffusion models
  for inverse problems through stochastic contraction.
\newblock In {\em Proceedings of the IEEE/CVF Conference on Computer Vision and
  Pattern Recognition (CVPR)}, pages 12413--12422, 2022.

\bibitem{chung2022score}
Hyungjin Chung and Jong~Chul Ye.
\newblock Score-based diffusion models for accelerated {MRI}.
\newblock {\em Medical Image Analysis}, 80:102479, 2022.

\bibitem{cohen2021regularization}
Regev Cohen, Michael Elad, and Peyman Milanfar.
\newblock Regularization by denoising via fixed-point projection (red-pro).
\newblock {\em SIAM Journal on Imaging Sciences}, 14(3):1374--1406, 2021.

\bibitem{dabov2007image}
Kostadin Dabov, Alessandro Foi, Vladimir Katkovnik, and Karen Egiazarian.
\newblock Image denoising by sparse 3-d transform-domain collaborative
  filtering.
\newblock {\em IEEE Transactions on Image Processing}, 16(8):2080--2095, 2007.

\bibitem{darestani2021accelerated}
Mohammad~Zalbagi Darestani and Reinhard Heckel.
\newblock Accelerated {MRI} with un-trained neural networks.
\newblock {\em IEEE Transactions on Computational Imaging}, 7:724--733, 2021.

\bibitem{daubechies2004ista}
Ingrid Daubechies, Michel Defrise, and Christine De~Mol.
\newblock An iterative thresholding algorithm for linear inverse problems with
  a sparsity constraint.
\newblock {\em Communications on Pure and Applied Mathematics},
  57(11):1413--1457, 2004.

\bibitem{reyes2016total}
J.C. {De Los Reyes}, C.-B. Schönlieb, and T.~Valkonen.
\newblock The structure of optimal parameters for image restoration problems.
\newblock {\em Journal of Mathematical Analysis and Applications},
  434(1):464--500, 2016.

\bibitem{dhar2018modeling}
Manik Dhar, Aditya Grover, and Stefano Ermon.
\newblock Modeling sparse deviations for compressed sensing using generative
  models.
\newblock In Jennifer Dy and Andreas Krause, editors, {\em Proceedings of the
  35th International Conference on Machine Learning}, volume~80 of {\em
  Proceedings of Machine Learning Research}, pages 1214--1223. PMLR, 2018.

\bibitem{dinh2015nice}
Laurent Dinh, David Krueger, and Yoshua Bengio.
\newblock {NICE:} non-linear independent components estimation.
\newblock In Yoshua Bengio and Yann LeCun, editors, {\em 3rd International
  Conference on Learning Representations, {ICLR} 2015, San Diego, CA, USA, May
  7-9, 2015, Workshop Track Proceedings}, 2015.

\bibitem{dinh2017density}
Laurent Dinh, Jascha Sohl{-}Dickstein, and Samy Bengio.
\newblock Density estimation using real {NVP}.
\newblock In {\em 5th International Conference on Learning Representations,
  {ICLR} 2017, Toulon, France, April 24-26, 2017, Conference Track
  Proceedings}. OpenReview.net, 2017.

\bibitem{dittmer2020regularization}
S{\"o}ren Dittmer, Tobias Kluth, Peter Maass, and Daniel Otero~Baguer.
\newblock Regularization by architecture: A deep prior approach for inverse
  problems.
\newblock {\em Journal of Mathematical Imaging and Vision}, 62:456--470, 2020.

\bibitem{dong2019denoising}
Weisheng Dong, Peiyao Wang, Wotao Yin, Guangming Shi, Fangfang Wu, and Xiaotong
  Lu.
\newblock Denoising prior driven deep neural network for image restoration.
\newblock {\em IEEE Transactions on Pattern Analysis and Machine Intelligence},
  41(10):2305--2318, 2019.

\bibitem{duff2021regularising}
Margaret Duff, Neill~DF Campbell, and Matthias~J Ehrhardt.
\newblock Regularising inverse problems with generative machine learning
  models.
\newblock {\em arXiv preprint arXiv:2107.11191}, 2021.

\bibitem{duff2023vaes}
Margaret Duff, Ivor Simpson, Matthias~J Ehrhardt, and Neill~DF Campbell.
\newblock {VAEs} with structured image covariance applied to compressed sensing
  {MRI}.
\newblock {\em Physics in Medicine \& Biology}, 68(16):165008, 2023.

\bibitem{durmus2019analysis}
Alain Durmus, Szymon Majewski, and B{\l}a{\.z}ej Miasojedow.
\newblock Analysis of langevin monte carlo via convex optimization.
\newblock {\em The Journal of Machine Learning Research}, 20(1):2666--2711,
  2019.

\bibitem{durmus19langevin}
Alain Durmus and {\'E}ric Moulines.
\newblock {High-dimensional Bayesian inference via the unadjusted Langevin
  algorithm}.
\newblock {\em Bernoulli}, 25(4A):2854 -- 2882, 2019.

\bibitem{durmus2022proximal}
Alain Durmus, {\'E}ric Moulines, and Marcelo Pereyra.
\newblock A proximal markov chain monte carlo method for bayesian inference in
  imaging inverse problems: When langevin meets moreau.
\newblock {\em SIAM Review}, 64(4):991--1028, 2022.

\bibitem{ebner2022plug}
Andrea Ebner and Markus Haltmeier.
\newblock Plug-and-play image reconstruction is a convergent regularization
  method.
\newblock {\em arXiv preprint arXiv:2212.06881}, 2022.

\bibitem{efron2011tweedie}
Bradley Efron.
\newblock Tweedie’s formula and selection bias.
\newblock {\em Journal of the American Statistical Association},
  106(496):1602--1614, 2011.
\newblock PMID: 22505788.

\bibitem{elad2006image}
Michael Elad and Michal Aharon.
\newblock Image denoising via sparse and redundant representations over learned
  dictionaries.
\newblock {\em IEEE Transactions on Image Processing}, 15(12):3736--3745, 2006.

\bibitem{engl1996regularization}
Heinz~Werner Engl and Martin Hanke.
\newblock {\em Regularization of inverse problems}, volume 375.
\newblock Springer Science \& Business Media, 1996.

\bibitem{erlacher2023joint}
Moritz Erlacher and Martin Zach.
\newblock Joint non-linear {MRI} inversion with diffusion priors.
\newblock {\em arXiv preprint arXiv:2310.14842}, 2023.

\bibitem{gao2023fast}
Guangyu Gao, Bo~Han, Zhenwu Fu, and Shanshan Tong.
\newblock A fast data-driven iteratively regularized method with convex penalty
  for solving ill-posed problems.
\newblock {\em SIAM Journal on Imaging Sciences}, 16(2):640--670, 2023.

\bibitem{gonzales2022solving}
Mario Gonz{\'a}lez, Andr{\'e}s Almansa, and Pauline Tan.
\newblock Solving inverse problems by joint posterior maximization with
  autoencoding prior.
\newblock {\em SIAM Journal on Imaging Sciences}, 15(2):822--859, 2022.

\bibitem{goodfellow2014generative}
Ian Goodfellow, Jean Pouget-Abadie, Mehdi Mirza, Bing Xu, David Warde-Farley,
  Sherjil Ozair, Aaron Courville, and Yoshua Bengio.
\newblock Generative adversarial nets.
\newblock In Z.~Ghahramani, M.~Welling, C.~Cortes, N.~Lawrence, and K.Q.
  Weinberger, editors, {\em Advances in Neural Information Processing Systems},
  volume~27. Curran Associates, Inc., 2014.

\bibitem{gregor2010learning}
Karol Gregor and Yann LeCun.
\newblock Learning fast approximations of sparse coding.
\newblock In {\em Proceedings of the 27th International Conference on
  International Conference on Machine Learning}, ICML'10, page 399–406,
  Madison, WI, USA, 2010. Omnipress.

\bibitem{guan2023magnetic}
Yu~Guan, Zongjiang Tu, Shanshan Wang, Yuhao Wang, Qiegen Liu, and Dong Liang.
\newblock Magnetic resonance imaging reconstruction using a deep energy-based
  model.
\newblock {\em NMR in Biomedicine}, 36(3):e4848, 2023.

\bibitem{gulrajani2017improved}
Ishaan Gulrajani, Faruk Ahmed, Martin Arjovsky, Vincent Dumoulin, and Aaron~C
  Courville.
\newblock Improved training of wasserstein {GANs}.
\newblock In I.~Guyon, U.~Von Luxburg, S.~Bengio, H.~Wallach, R.~Fergus,
  S.~Vishwanathan, and R.~Garnett, editors, {\em Advances in Neural Information
  Processing Systems}, volume~30. Curran Associates, Inc., 2017.

\bibitem{guo2019agem}
Bichuan Guo, Yuxing Han, and Jiangtao Wen.
\newblock Agem: Solving linear inverse problems via deep priors and sampling.
\newblock In H.~Wallach, H.~Larochelle, A.~Beygelzimer, F.~d\textquotesingle
  Alch\'{e}-Buc, E.~Fox, and R.~Garnett, editors, {\em Advances in Neural
  Information Processing Systems}, volume~32. Curran Associates, Inc., 2019.

\bibitem{gurrola2021residual}
Javier Gurrola-Ramos, Oscar Dalmau, and Teresa~E. Alarcón.
\newblock A residual dense u-net neural network for image denoising.
\newblock {\em IEEE Access}, 9:31742--31754, 2021.

\bibitem{habring2022generative}
Andreas Habring and Martin Holler.
\newblock A generative variational model for inverse problems in imaging.
\newblock {\em SIAM Journal on Mathematics of Data Science}, 4(1):306--335,
  2022.

\bibitem{habring2023note}
Andreas Habring and Martin Holler.
\newblock A note on the regularity of images generated by convolutional neural
  networks.
\newblock {\em SIAM Journal on Mathematics of Data Science}, 5(3):670--692,
  2023.

\bibitem{habring2023subgradient}
Andreas Habring, Martin Holler, and Thomas Pock.
\newblock Subgradient langevin methods for sampling from non-smooth potentials.
\newblock {\em arXiv preprint arXiv:2308.01417}, 2023.

\bibitem{hagemann2023multilevel}
Paul Hagemann, Lars Ruthotto, Gabriele Steidl, and Nicole~Tianjiao Yang.
\newblock Multilevel diffusion: Infinite dimensional score-based diffusion
  models for image generation.
\newblock {\em arXiv preprint arXiv:2303.04772}, 2023.

\bibitem{haltmeier2023data}
Markus Haltmeier, Richard Kowar, and Markus Tiefentaler.
\newblock Data-driven {Morozov} regularization of inverse problems.
\newblock {\em arXiv preprint arXiv:2310.14290}, 2023.

\bibitem{hammernik2018learning}
Kerstin Hammernik, Teresa Klatzer, Erich Kobler, Michael~P. Recht, Daniel~K.
  Sodickson, Thomas Pock, and Florian Knoll.
\newblock Learning a variational network for reconstruction of accelerated
  {MRI} data.
\newblock {\em Magnetic Resonance in Medicine}, 79(6):3055--3071, 2018.

\bibitem{hand2018phase}
Paul Hand, Oscar Leong, and Vlad Voroninski.
\newblock Phase retrieval under a generative prior.
\newblock In S.~Bengio, H.~Wallach, H.~Larochelle, K.~Grauman, N.~Cesa-Bianchi,
  and R.~Garnett, editors, {\em Advances in Neural Information Processing
  Systems}, volume~31. Curran Associates, Inc., 2018.

\bibitem{hasannasab2020parseval}
Marzieh Hasannasab, Johannes Hertrich, Sebastian Neumayer, Gerlind Plonka,
  Simon Setzer, and Gabriele Steidl.
\newblock Parseval proximal neural networks.
\newblock {\em Journal of Fourier Analysis and Applications}, 26:1--31, 2020.

\bibitem{he2019optimizing}
Ji~He, Yan Yang, Yongbo Wang, Dong Zeng, Zhaoying Bian, Hao Zhang, Jian Sun,
  Zongben Xu, and Jianhua Ma.
\newblock Optimizing a parameterized plug-and-play admm for iterative low-dose
  ct reconstruction.
\newblock {\em IEEE Transactions on Medical Imaging}, 38(2):371--382, 2019.

\bibitem{heckel2019regularizing}
Reinhard Heckel.
\newblock Regularizing linear inverse problems with convolutional neural
  networks.
\newblock {\em arXiv preprint arXiv:1907.03100}, 2019.

\bibitem{heckel2018deep}
Reinhard Heckel and Paul Hand.
\newblock Deep decoder: Concise image representations from untrained
  non-convolutional networks.
\newblock {\em arXiv preprint arXiv:1810.03982}, 2018.

\bibitem{heckel2019denoising}
Reinhard Heckel and Mahdi Soltanolkotabi.
\newblock Denoising and regularization via exploiting the structural bias of
  convolutional generators.
\newblock {\em arXiv preprint arXiv:1910.14634}, 2019.

\bibitem{heckel2020compressive}
Reinhard Heckel and Mahdi Soltanolkotabi.
\newblock Compressive sensing with un-trained neural networks: Gradient descent
  finds a smooth approximation.
\newblock In {\em International Conference on Machine Learning}, pages
  4149--4158. PMLR, 2020.

\bibitem{hendriksen2020noise}
Allard~Adriaan Hendriksen, Daniël~Maria Pelt, and K.~Joost Batenburg.
\newblock Noise2inverse: Self-supervised deep convolutional denoising for
  tomography.
\newblock {\em IEEE Transactions on Computational Imaging}, 6:1320--1335, 2020.

\bibitem{hinton2002training}
Geoffrey~E. Hinton.
\newblock Training products of experts by minimizing contrastive divergence.
\newblock {\em Neural Computation}, 14(8):1771--1800, 2002.

\bibitem{ho2020denoising}
Jonathan Ho, Ajay Jain, and Pieter Abbeel.
\newblock Denoising diffusion probabilistic models.
\newblock In H.~Larochelle, M.~Ranzato, R.~Hadsell, M.F. Balcan, and H.~Lin,
  editors, {\em Advances in Neural Information Processing Systems}, volume~33,
  pages 6840--6851. Curran Associates, Inc., 2020.

\bibitem{holden2022bayesian}
Matthew Holden, Marcelo Pereyra, and Konstantinos~C Zygalakis.
\newblock Bayesian imaging with data-driven priors encoded by neural networks.
\newblock {\em SIAM Journal on Imaging Sciences}, 15(2):892--924, 2022.

\bibitem{hong2019acceleration}
Tao Hong, Yaniv Romano, and Michael Elad.
\newblock Acceleration of red via vector extrapolation.
\newblock {\em Journal of Visual Communication and Image Representation},
  63:102575, 2019.

\bibitem{chang2018deep}
Chang~Min Hyun, Hwa~Pyung Kim, Sung~Min Lee, Sungchul Lee, and Jin~Keun Seo.
\newblock Deep learning for undersampled {MRI} reconstruction.
\newblock {\em Physics in Medicine \& Biology}, 63(13):135007, 2018.

\bibitem{hyvarinen2005estimation}
Aapo Hyv{{\"a}}rinen.
\newblock Estimation of non-normalized statistical models by score matching.
\newblock {\em Journal of Machine Learning Research}, 6(24):695--709, 2005.

\bibitem{isola2017image}
Phillip Isola, Jun-Yan Zhu, Tinghui Zhou, and Alexei~A. Efros.
\newblock Image-to-image translation with conditional adversarial networks.
\newblock In {\em Proceedings of the IEEE Conference on Computer Vision and
  Pattern Recognition (CVPR)}, 2017.

\bibitem{Ivanov65_ivanov_regularization}
Valentin~Konstantinovich Ivanov.
\newblock On linear problems which are not well-posed.
\newblock {\em Doklady Akademii Nauk SSSR}, 145(2):270--272, 1962.

\bibitem{jagatap2019algorithmic}
Gauri Jagatap and Chinmay Hegde.
\newblock Algorithmic guarantees for inverse imaging with untrained network
  priors.
\newblock {\em Advances in neural information processing systems}, 32, 2019.

\bibitem{jain2008natural}
Viren Jain and Sebastian Seung.
\newblock Natural image denoising with convolutional networks.
\newblock In D.~Koller, D.~Schuurmans, Y.~Bengio, and L.~Bottou, editors, {\em
  Advances in Neural Information Processing Systems}, volume~21. Curran
  Associates, Inc., 2008.

\bibitem{jin2017deep}
Kyong~Hwan Jin, Michael~T. McCann, Emmanuel Froustey, and Michael Unser.
\newblock Deep convolutional neural network for inverse problems in imaging.
\newblock {\em IEEE Transactions on Image Processing}, 26(9):4509--4522, 2017.

\bibitem{kabri2022convergent}
Samira Kabri, Alexander Auras, Danilo Riccio, Hartmut Bauermeister, Martin
  Benning, Michael Moeller, and Martin Burger.
\newblock Convergent data-driven regularizations for ct reconstruction.
\newblock {\em arXiv preprint arXiv:2212.07786}, 2022.

\bibitem{kamilov2017plug}
Ulugbek~S. Kamilov, Hassan Mansour, and Brendt Wohlberg.
\newblock A plug-and-play priors approach for solving nonlinear imaging inverse
  problems.
\newblock {\em IEEE Signal Processing Letters}, 24(12):1872--1876, 2017.

\bibitem{kang2017deep}
Eunhee Kang, Junhong Min, and Jong~Chul Ye.
\newblock A deep convolutional neural network using directional wavelets for
  low-dose x-ray ct reconstruction.
\newblock {\em Medical Physics}, 44(10):e360--e375, 2017.

\bibitem{karami2019invertible}
Mahdi Karami, Dale Schuurmans, Jascha Sohl-Dickstein, Laurent Dinh, and Daniel
  Duckworth.
\newblock Invertible convolutional flow.
\newblock In H.~Wallach, H.~Larochelle, A.~Beygelzimer, F.~d\textquotesingle
  Alch\'{e}-Buc, E.~Fox, and R.~Garnett, editors, {\em Advances in Neural
  Information Processing Systems}, volume~32. Curran Associates, Inc., 2019.

\bibitem{kawar2022denoising}
Bahjat Kawar, Michael Elad, Stefano Ermon, and Jiaming Song.
\newblock Denoising diffusion restoration models.
\newblock In S.~Koyejo, S.~Mohamed, A.~Agarwal, D.~Belgrave, K.~Cho, and A.~Oh,
  editors, {\em Advances in Neural Information Processing Systems}, volume~35,
  pages 23593--23606. Curran Associates, Inc., 2022.

\bibitem{kawar2021snips}
Bahjat Kawar, Gregory Vaksman, and Michael Elad.
\newblock Snips: Solving noisy inverse problems stochastically.
\newblock In M.~Ranzato, A.~Beygelzimer, Y.~Dauphin, P.S. Liang, and J.~Wortman
  Vaughan, editors, {\em Advances in Neural Information Processing Systems},
  volume~34, pages 21757--21769. Curran Associates, Inc., 2021.

\bibitem{kervrann2006optimal}
C.~Kervrann and J.~Boulanger.
\newblock Optimal spatial adaptation for patch-based image denoising.
\newblock {\em IEEE Transactions on Image Processing}, 15(10):2866--2878, 2006.

\bibitem{kingma2013auto}
Diederik~P. Kingma and Max Welling.
\newblock Auto-encoding variational bayes.
\newblock In Yoshua Bengio and Yann LeCun, editors, {\em 2nd International
  Conference on Learning Representations, {ICLR} 2014, Banff, AB, Canada, April
  14-16, 2014, Conference Track Proceedings}, 2014.

\bibitem{kingma2018glow}
Durk~P Kingma and Prafulla Dhariwal.
\newblock Glow: Generative flow with invertible 1x1 convolutions.
\newblock In S.~Bengio, H.~Wallach, H.~Larochelle, K.~Grauman, N.~Cesa-Bianchi,
  and R.~Garnett, editors, {\em Advances in Neural Information Processing
  Systems}, volume~31. Curran Associates, Inc., 2018.

\bibitem{knoll2020fastmri}
Florian Knoll, Jure Zbontar, Anuroop Sriram, Matthew~J Muckley, Mary Bruno,
  Aaron Defazio, Marc Parente, Krzysztof~J Geras, Joe Katsnelson, Hersh
  Chandarana, et~al.
\newblock {fastMRI}: A publicly available raw k-space and {DICOM} dataset of
  knee images for accelerated {MR} image reconstruction using machine learning.
\newblock {\em Radiology: Artificial Intelligence}, 2(1):e190007, 2020.

\bibitem{kobler2020total}
Erich Kobler, Alexander Effland, Karl Kunisch, and Thomas Pock.
\newblock Total deep variation for linear inverse problems.
\newblock In {\em Proceedings of the IEEE/CVF Conference on computer vision and
  pattern recognition}, pages 7549--7558, 2020.

\bibitem{kobler2022total}
Erich Kobler, Alexander Effland, Karl Kunisch, and Thomas Pock.
\newblock Total deep variation: A stable regularization method for inverse
  problems.
\newblock {\em IEEE Transactions on Pattern Analysis and Machine Intelligence},
  44(12):9163--9180, 2022.

\bibitem{korolev2022two}
Yury Korolev.
\newblock Two-layer neural networks with values in a banach space.
\newblock {\em SIAM Journal on Mathematical Analysis}, 54(6):6358--6389, 2022.

\bibitem{kunisch2013bilevel}
Karl Kunisch and Thomas Pock.
\newblock A bilevel optimization approach for parameter learning in variational
  models.
\newblock {\em SIAM Journal on Imaging Sciences}, 6(2):938--983, 2013.

\bibitem{lanthaler2022error}
Samuel Lanthaler, Siddhartha Mishra, and George~E Karniadakis.
\newblock Error estimates for deeponets: A deep learning framework in infinite
  dimensions.
\newblock {\em Transactions of Mathematics and Its Applications}, 6(1):tnac001,
  2022.

\bibitem{laumont2022bayesian}
R\'{e}mi Laumont, Valentin~De Bortoli, Andr\'{e}s Almansa, Julie Delon, Alain
  Durmus, and Marcelo Pereyra.
\newblock Bayesian imaging using plug \& play priors: When langevin meets
  tweedie.
\newblock {\em SIAM Journal on Imaging Sciences}, 15(2):701--737, 2022.

\bibitem{lehtinen2018noise}
Jaakko Lehtinen, Jacob Munkberg, Jon Hasselgren, Samuli Laine, Tero Karras,
  Miika Aittala, and Timo Aila.
\newblock {N}oise2{N}oise: Learning image restoration without clean data.
\newblock In Jennifer Dy and Andreas Krause, editors, {\em Proceedings of the
  35th International Conference on Machine Learning}, volume~80 of {\em
  Proceedings of Machine Learning Research}, pages 2965--2974. PMLR, 2018.

\bibitem{li2020nett}
Housen Li, Johannes Schwab, Stephan Antholzer, and Markus Haltmeier.
\newblock Nett: solving inverse problems with deep neural networks.
\newblock {\em Inverse Problems}, 36(6):065005, 2020.

\bibitem{li2021review}
Y~Li, Bruno Sixou, and F~Peyrin.
\newblock A review of the deep learning methods for medical images super
  resolution problems.
\newblock {\em IRBM}, 42(2):120--133, 2021.

\bibitem{liu2019image}
Jiaming Liu, Yu~Sun, Xiaojian Xu, and Ulugbek~S Kamilov.
\newblock Image restoration using total variation regularized deep image prior.
\newblock In {\em ICASSP 2019-2019 IEEE International Conference on Acoustics,
  Speech and Signal Processing (ICASSP)}, pages 7715--7719. Ieee, 2019.

\bibitem{Liu_2023_ICCV}
Yilin Liu, Jiang Li, Yunkui Pang, Dong Nie, and Pew-Thian Yap.
\newblock The devil is in the upsampling: Architectural decisions made simpler
  for denoising with deep image prior.
\newblock In {\em Proceedings of the IEEE/CVF International Conference on
  Computer Vision (ICCV)}, pages 12408--12417, 2023.

\bibitem{lucas2018using}
Alice Lucas, Michael Iliadis, Rafael Molina, and Aggelos~K Katsaggelos.
\newblock Using deep neural networks for inverse problems in imaging: beyond
  analytical methods.
\newblock {\em IEEE Signal Processing Magazine}, 35(1):20--36, 2018.

\bibitem{lunz2018adversarial}
Sebastian Lunz, Ozan \"{O}ktem, and Carola-Bibiane Sch\"{o}nlieb.
\newblock Adversarial regularizers in inverse problems.
\newblock In S.~Bengio, H.~Wallach, H.~Larochelle, K.~Grauman, N.~Cesa-Bianchi,
  and R.~Garnett, editors, {\em Advances in Neural Information Processing
  Systems}, volume~31. Curran Associates, Inc., 2018.

\bibitem{luo2023bayesian}
Guanxiong Luo, Moritz Blumenthal, Martin Heide, and Martin Uecker.
\newblock Bayesian {MRI} reconstruction with joint uncertainty estimation using
  diffusion models.
\newblock {\em Magnetic Resonance in Medicine}, 90(1):295--311, 2023.

\bibitem{mairal2009nonlocal}
Julien Mairal, Francis Bach, Jean Ponce, Guillermo Sapiro, and Andrew
  Zisserman.
\newblock Non-local sparse models for image restoration.
\newblock In {\em 2009 IEEE 12th International Conference on Computer Vision},
  pages 2272--2279, 2009.

\bibitem{mataev2019deepred}
Gary Mataev, Peyman Milanfar, and Michael Elad.
\newblock Deepred: Deep image prior powered by red.
\newblock In {\em Proceedings of the IEEE/CVF International Conference on
  Computer Vision Workshops}, pages 0--0, 2019.

\bibitem{mccann2017convolutional}
Michael~T McCann, Kyong~Hwan Jin, and Michael Unser.
\newblock Convolutional neural networks for inverse problems in imaging: A
  review.
\newblock {\em IEEE Signal Processing Magazine}, 34(6):85--95, 2017.

\bibitem{meinhardt2017learning}
Tim Meinhardt, Michael Moller, Caner Hazirbas, and Daniel Cremers.
\newblock Learning proximal operators: Using denoising networks for
  regularizing inverse imaging problems.
\newblock In {\em Proceedings of the IEEE International Conference on Computer
  Vision (ICCV)}, 2017.

\bibitem{metzler2018prdeep}
Christopher Metzler, Phillip Schniter, Ashok Veeraraghavan, and Richard
  Baraniuk.
\newblock pr{D}eep: Robust phase retrieval with a flexible deep network.
\newblock In Jennifer Dy and Andreas Krause, editors, {\em Proceedings of the
  35th International Conference on Machine Learning}, volume~80 of {\em
  Proceedings of Machine Learning Research}, pages 3501--3510. PMLR, 2018.

\bibitem{moreau1965proximite}
Jean-Jacques Moreau.
\newblock Proximit{\'e} et dualit{\'e} dans un espace hilbertien.
\newblock {\em Bulletin de la Soci{\'e}t{\'e} math{\'e}matique de France},
  93:273--299, 1965.

\bibitem{mosser2020stochastic}
Lukas Mosser, Olivier Dubrule, and Martin~J. Blunt.
\newblock Stochastic seismic waveform inversion using generative adversarial
  networks as a geological prior.
\newblock {\em Mathematical Geosciences}, 52(1):53--79, 2020.

\bibitem{mukherjee2020learned}
Subhadip Mukherjee, S{\"o}ren Dittmer, Zakhar Shumaylov, Sebastian Lunz, Ozan
  {\"O}ktem, and Carola-Bibiane Sch{\"o}nlieb.
\newblock Learned convex regularizers for inverse problems.
\newblock {\em arXiv preprint arXiv:2008.02839}, 2020.

\bibitem{mukherjee2023learned}
Subhadip Mukherjee, Andreas Hauptmann, Ozan {\"O}ktem, Marcelo Pereyra, and
  Carola-Bibiane Sch{\"o}nlieb.
\newblock Learned reconstruction methods with convergence guarantees: a survey
  of concepts and applications.
\newblock {\em IEEE Signal Processing Magazine}, 40(1):164--182, 2023.

\bibitem{narnhofer2022bayesian}
Dominik Narnhofer, Alexander Effland, Erich Kobler, Kerstin Hammernik, Florian
  Knoll, and Thomas Pock.
\newblock Bayesian uncertainty estimation of learned variational {MRI}
  reconstruction.
\newblock {\em IEEE Transactions on Medical Imaging}, 41(2):279--291, 2022.

\bibitem{narnhofer2022posterior}
Dominik Narnhofer, Andreas Habring, Martin Holler, and Thomas Pock.
\newblock Posterior-variance-based error quantification for inverse problems in
  imaging.
\newblock {\em arXiv preprint arXiv:2212.12499}, 2022.

\bibitem{obmann2021augmented}
Daniel Obmann, Linh Nguyen, Johannes Schwab, and Markus Haltmeier.
\newblock Augmented nett regularization of inverse problems.
\newblock {\em Journal of Physics Communications}, 5(10):105002, 2021.

\bibitem{obmann2020deep}
Daniel Obmann, Johannes Schwab, and Markus Haltmeier.
\newblock Deep synthesis network for regularizing inverse problems.
\newblock {\em Inverse Problems}, 37(1):015005, 2020.

\bibitem{ongie2020deep}
Gregory Ongie, Ajil Jalal, Christopher~A Metzler, Richard~G Baraniuk,
  Alexandros~G Dimakis, and Rebecca Willett.
\newblock Deep learning techniques for inverse problems in imaging.
\newblock {\em IEEE Journal on Selected Areas in Information Theory},
  1(1):39--56, 2020.

\bibitem{pan2022exploiting}
Xingang Pan, Xiaohang Zhan, Bo~Dai, Dahua Lin, Chen~Change Loy, and Ping Luo.
\newblock Exploiting deep generative prior for versatile image restoration and
  manipulation.
\newblock {\em IEEE Transactions on Pattern Analysis and Machine Intelligence},
  44(11):7474--7489, 2022.

\bibitem{pelt2018improving}
Daniël~M. Pelt, Kees~Joost Batenburg, and James~A. Sethian.
\newblock Improving tomographic reconstruction from limited data using
  mixed-scale dense convolutional neural networks.
\newblock {\em Journal of Imaging}, 4(11), 2018.

\bibitem{putzky2017recurrent}
Patrick Putzky and Max Welling.
\newblock Recurrent inference machines for solving inverse problems.
\newblock {\em arXiv preprint arXiv:1706.04008}, 2017.

\bibitem{qayyum2022untrained}
Adnan Qayyum, Inaam Ilahi, Fahad Shamshad, Farid Boussaid, Mohammed Bennamoun,
  and Junaid Qadir.
\newblock Untrained neural network priors for inverse imaging problems: A
  survey.
\newblock {\em IEEE Transactions on Pattern Analysis and Machine Intelligence},
  2022.

\bibitem{quan2018compressed}
Tran~Minh Quan, Thanh Nguyen-Duc, and Won-Ki Jeong.
\newblock Compressed sensing {MRI} reconstruction using a generative
  adversarial network with a cyclic loss.
\newblock {\em IEEE Transactions on Medical Imaging}, 37(6):1488--1497, 2018.

\bibitem{radford2016dcgan}
Alec Radford, Luke Metz, and Soumith Chintala.
\newblock Unsupervised representation learning with deep convolutional
  generative adversarial networks.
\newblock In Yoshua Bengio and Yann LeCun, editors, {\em 4th International
  Conference on Learning Representations, {ICLR} 2016, San Juan, Puerto Rico,
  May 2-4, 2016, Conference Track Proceedings}, 2016.

\bibitem{raj2019ganbased}
Ankit Raj, Yuqi Li, and Yoram Bresler.
\newblock Gan-based projector for faster recovery with convergence guarantees
  in linear inverse problems.
\newblock In {\em Proceedings of the IEEE/CVF International Conference on
  Computer Vision (ICCV)}, 2019.

\bibitem{ramzi2020denoising}
Zaccharie Ramzi, Benjamin Remy, Francois Lanusse, Jean-Luc Starck, and Philippe
  Ciuciu.
\newblock Denoising score-matching for uncertainty quantification in inverse
  problems.
\newblock {\em arXiv preprint arXiv:2011.08698}, 2020.

\bibitem{reehorst2019red}
Edward~T. Reehorst and Philip Schniter.
\newblock Regularization by denoising: Clarifications and new interpretations.
\newblock {\em IEEE Transactions on Computational Imaging}, 5(1):52--67, 2019.

\bibitem{rezende2015variational}
Danilo Rezende and Shakir Mohamed.
\newblock Variational inference with normalizing flows.
\newblock In Francis Bach and David Blei, editors, {\em Proceedings of the 32nd
  International Conference on Machine Learning}, volume~37 of {\em Proceedings
  of Machine Learning Research}, pages 1530--1538, Lille, France, 2015. PMLR.

\bibitem{chang2017one}
J.~H. Rick~Chang, Chun-Liang Li, Barnabas Poczos, B.~V.~K. Vijaya~Kumar, and
  Aswin~C. Sankaranarayanan.
\newblock One network to solve them all -- solving linear inverse problems
  using deep projection models.
\newblock In {\em Proceedings of the IEEE International Conference on Computer
  Vision (ICCV)}, 2017.

\bibitem{rizzuti2020param}
Gabrio Rizzuti, Ali Siahkoohi, Philipp~A. Witte, and Felix~J. Herrmann.
\newblock Parameterizing uncertainty by deep invertible networks: An
  application to reservoir characterization.
\newblock In {\em SEG International Exposition and Annual Meeting}, page
  D031S057R006, 2020.

\bibitem{romano2017red}
Yaniv Romano, Michael Elad, and Peyman Milanfar.
\newblock The little engine that could: Regularization by denoising (red).
\newblock {\em SIAM Journal on Imaging Sciences}, 10(4):1804--1844, 2017.

\bibitem{romano2019conformalized}
Yaniv Romano, Evan Patterson, and Emmanuel Candes.
\newblock Conformalized quantile regression.
\newblock {\em Advances in neural information processing systems}, 32, 2019.

\bibitem{ronneberger2015u}
Olaf Ronneberger, Philipp Fischer, and Thomas Brox.
\newblock U-net: Convolutional networks for biomedical image segmentation.
\newblock In {\em International Conference on Medical image computing and
  computer-assisted intervention}, pages 234--241. Springer, 2015.

\bibitem{roth2005fields}
Stefan Roth and Michael~J Black.
\newblock Fields of experts: A framework for learning image priors.
\newblock In {\em 2005 IEEE Computer Society Conference on Computer Vision and
  Pattern Recognition (CVPR'05)}, volume~2, pages 860--867. IEEE, 2005.

\bibitem{ro09}
Stefan Roth and Michael~J Black.
\newblock Fields of experts.
\newblock {\em International Journal of Computer Vision}, 82(2):205--229, 2009.

\bibitem{rudin1992nonlinear}
Leonid~I Rudin, Stanley Osher, and Emad Fatemi.
\newblock Nonlinear total variation based noise removal algorithms.
\newblock {\em Physica D: nonlinear phenomena}, 60(1-4):259--268, 1992.

\bibitem{rudin1992tv_mh}
Leonid~I. Rudin, Stanley Osher, and Emad Fatemi.
\newblock Nonlinear total variation based noise removal algorithms.
\newblock {\em Physica D}, 60(1--4):259--268, 1992.

\bibitem{ryu2019plug}
Ernest Ryu, Jialin Liu, Sicheng Wang, Xiaohan Chen, Zhangyang Wang, and Wotao
  Yin.
\newblock Plug-and-play methods provably converge with properly trained
  denoisers.
\newblock In Kamalika Chaudhuri and Ruslan Salakhutdinov, editors, {\em
  Proceedings of the 36th International Conference on Machine Learning},
  volume~97 of {\em Proceedings of Machine Learning Research}, pages
  5546--5557. PMLR, 2019.

\bibitem{scarlett2022theoretical}
Jonathan Scarlett, Reinhard Heckel, Miguel~RD Rodrigues, Paul Hand, and
  Yonina~C Eldar.
\newblock Theoretical perspectives on deep learning methods in inverse
  problems.
\newblock {\em IEEE journal on selected areas in information theory},
  3(3):433--453, 2022.

\bibitem{scherzer2008variational}
Otmar Scherzer, Markus Grasmair, Harald Grossauer, Markus Haltmeier, and Frank
  Lenzen.
\newblock {\em Variational Methods in Imaging}, volume 167.
\newblock Springer Science \& Business Media, 2008.

\bibitem{scherzer2023gauss}
Otmar Scherzer, Bernd Hofmann, and Zuhair Nashed.
\newblock Gauss--newton method for solving linear inverse problems with neural
  network coders.
\newblock {\em Sampling Theory, Signal Processing, and Data Analysis},
  21(2):25, 2023.

\bibitem{schlemper2018deep}
Jo~Schlemper, Jose Caballero, Joseph~V. Hajnal, Anthony~N. Price, and Daniel
  Rueckert.
\newblock A deep cascade of convolutional neural networks for dynamic mr image
  reconstruction.
\newblock {\em IEEE Transactions on Medical Imaging}, 37(2):491--503, 2018.

\bibitem{schlemper2018bayesian}
Jo~Schlemper, Daniel~C Castro, Wenjia Bai, Chen Qin, Ozan Oktay, Jinming Duan,
  Anthony~N Price, Jo~Hajnal, and Daniel Rueckert.
\newblock Bayesian deep learning for accelerated mr image reconstruction.
\newblock In {\em Machine Learning for Medical Image Reconstruction: First
  International Workshop, MLMIR 2018, Held in Conjunction with MICCAI 2018,
  Granada, Spain, September 16, 2018, Proceedings 1}, pages 64--71. Springer,
  2018.

\bibitem{schwab2019deep}
Johannes Schwab, Stephan Antholzer, and Markus Haltmeier.
\newblock Deep null space learning for inverse problems: convergence analysis
  and rates.
\newblock {\em Inverse Problems}, 35(2):025008, 2019.

\bibitem{schwab2020big}
Johannes Schwab, Stephan Antholzer, and Markus Haltmeier.
\newblock Big in japan: Regularizing networks for solving inverse problems.
\newblock {\em Journal of mathematical imaging and vision}, 62(3):445--455,
  2020.

\bibitem{senouf2019self}
Ortal Senouf, Sanketh Vedula, Tomer Weiss, Alex Bronstein, Oleg Michailovich,
  and Michael Zibulevsky.
\newblock Self-supervised learning of inverse problem solvers in medical
  imaging.
\newblock In Qian Wang, Fausto Milletari, Hien~V. Nguyen, Shadi Albarqouni,
  M.~Jorge Cardoso, Nicola Rieke, Ziyue Xu, Konstantinos Kamnitsas, Vishal
  Patel, Badri Roysam, Steve Jiang, Kevin Zhou, Khoa Luu, and Ngan Le, editors,
  {\em Domain Adaptation and Representation Transfer and Medical Image Learning
  with Less Labels and Imperfect Data}, pages 111--119, Cham, 2019. Springer
  International Publishing.

\bibitem{seo2019learning}
Jin~Keun Seo, Kang~Cheol Kim, Ariungerel Jargal, Kyounghun Lee, and Bastian
  Harrach.
\newblock A learning-based method for solving ill-posed nonlinear inverse
  problems: A simulation study of lung eit.
\newblock {\em SIAM Journal on Imaging Sciences}, 12(3):1275--1295, 2019.

\bibitem{shah2018solving}
Viraj Shah and Chinmay Hegde.
\newblock Solving linear inverse problems using gan priors: An algorithm with
  provable guarantees.
\newblock In {\em 2018 IEEE International Conference on Acoustics, Speech and
  Signal Processing (ICASSP)}, pages 4609--4613, 2018.

\bibitem{shlezinger2023model}
Nir Shlezinger, Jay Whang, Yonina~C Eldar, and Alexandros~G Dimakis.
\newblock Model-based deep learning.
\newblock {\em Proceedings of the IEEE}, 2023.

\bibitem{siahkoohi2020faster}
Ali Siahkoohi, Gabrio Rizzuti, Philipp~A Witte, and Felix~J Herrmann.
\newblock Faster uncertainty quantification for inverse problems with
  conditional normalizing flows.
\newblock {\em arXiv preprint arXiv:2007.07985}, 2020.

\bibitem{sohl2015deep}
Jascha Sohl-Dickstein, Eric Weiss, Niru Maheswaranathan, and Surya Ganguli.
\newblock Deep unsupervised learning using nonequilibrium thermodynamics.
\newblock In Francis Bach and David Blei, editors, {\em Proceedings of the 32nd
  International Conference on Machine Learning}, volume~37 of {\em Proceedings
  of Machine Learning Research}, pages 2256--2265, Lille, France, 2015. PMLR.

\bibitem{song2019generative}
Yang Song and Stefano Ermon.
\newblock Generative modeling by estimating gradients of the data distribution.
\newblock In H.~Wallach, H.~Larochelle, A.~Beygelzimer, F.~d\textquotesingle
  Alch\'{e}-Buc, E.~Fox, and R.~Garnett, editors, {\em Advances in Neural
  Information Processing Systems}, volume~32. Curran Associates, Inc., 2019.

\bibitem{song2020sliced}
Yang Song, Sahaj Garg, Jiaxin Shi, and Stefano Ermon.
\newblock Sliced score matching: A scalable approach to density and score
  estimation.
\newblock In Ryan~P. Adams and Vibhav Gogate, editors, {\em Proceedings of The
  35th Uncertainty in Artificial Intelligence Conference}, volume 115 of {\em
  Proceedings of Machine Learning Research}, pages 574--584. PMLR, 2020.

\bibitem{song2021solving}
Yang Song, Liyue Shen, Lei Xing, and Stefano Ermon.
\newblock Solving inverse problems in medical imaging with score-based
  generative models.
\newblock {\em arXiv preprint arXiv:2111.08005}, 2021.

\bibitem{song2020score}
Yang Song, Jascha Sohl-Dickstein, Diederik~P Kingma, Abhishek Kumar, Stefano
  Ermon, and Ben Poole.
\newblock Score-based generative modeling through stochastic differential
  equations.
\newblock {\em arXiv preprint arXiv:2011.13456}, 2020.

\bibitem{sreehari2016plug}
Suhas Sreehari, S.~V. Venkatakrishnan, Brendt Wohlberg, Gregery~T. Buzzard,
  Lawrence~F. Drummy, Jeffrey~P. Simmons, and Charles~A. Bouman.
\newblock Plug-and-play priors for bright field electron tomography and sparse
  interpolation.
\newblock {\em IEEE Transactions on Computational Imaging}, 2(4):408--423,
  2016.

\bibitem{stuart2010inverse}
Andrew~M Stuart.
\newblock Inverse problems: a bayesian perspective.
\newblock {\em Acta numerica}, 19:451--559, 2010.

\bibitem{sun2021deep}
He~Sun and Katherine~L. Bouman.
\newblock Deep probabilistic imaging: Uncertainty quantification and
  multi-modal solution characterization for computational imaging.
\newblock {\em Proceedings of the AAAI Conference on Artificial Intelligence},
  35(3):2628--2637, 2021.

\bibitem{sun2019block}
Yu~Sun, Jiaming Liu, and Ulugbek Kamilov.
\newblock Block coordinate regularization by denoising.
\newblock In H.~Wallach, H.~Larochelle, A.~Beygelzimer, F.~d\textquotesingle
  Alch\'{e}-Buc, E.~Fox, and R.~Garnett, editors, {\em Advances in Neural
  Information Processing Systems}, volume~32. Curran Associates, Inc., 2019.

\bibitem{sun2019online}
Yu~Sun, Brendt Wohlberg, and Ulugbek~S. Kamilov.
\newblock An online plug-and-play algorithm for regularized image
  reconstruction.
\newblock {\em IEEE Transactions on Computational Imaging}, 5(3):395--408,
  2019.

\bibitem{tasdizen2009principal}
Tolga Tasdizen.
\newblock Principal neighborhood dictionaries for nonlocal means image
  denoising.
\newblock {\em IEEE Transactions on Image Processing}, 18(12):2649--2660, 2009.

\bibitem{tenorio2017introduction}
Luis Tenorio.
\newblock {\em An introduction to data analysis and uncertainty quantification
  for inverse problems}.
\newblock SIAM, 2017.

\bibitem{tian2019enhanced}
Chunwei Tian, Yong Xu, Lunke Fei, Junqian Wang, Jie Wen, and Nan Luo.
\newblock Enhanced cnn for image denoising.
\newblock {\em CAAI Transactions on Intelligence Technology}, 4(1):17--23,
  2019.

\bibitem{tirer2019image}
Tom Tirer and Raja Giryes.
\newblock Image restoration by iterative denoising and backward projections.
\newblock {\em IEEE Transactions on Image Processing}, 28(3):1220--1234, 2019.

\bibitem{zongjiang2023kspace}
Zongjiang Tu, Chen Jiang, Yu~Guan, Jijun Liu, and Qiegen Liu.
\newblock K-space and image domain collaborative energy-based model for
  parallel {MRI} reconstruction.
\newblock {\em Magnetic Resonance Imaging}, 99:110--122, 2023.

\bibitem{ulyanov2018deep}
Dmitry Ulyanov, Andrea Vedaldi, and Victor Lempitsky.
\newblock Deep image prior.
\newblock In {\em Proceedings of the IEEE conference on computer vision and
  pattern recognition}, pages 9446--9454, 2018.

\bibitem{van2018compressed}
Dave Van~Veen, Ajil Jalal, Mahdi Soltanolkotabi, Eric Price, Sriram Vishwanath,
  and Alexandros~G Dimakis.
\newblock Compressed sensing with deep image prior and learned regularization.
\newblock {\em arXiv preprint arXiv:1806.06438}, 2018.

\bibitem{venkatakrishnan2013plug}
Singanallur~V. Venkatakrishnan, Charles~A. Bouman, and Brendt Wohlberg.
\newblock Plug-and-play priors for model based reconstruction.
\newblock In {\em 2013 IEEE Global Conference on Signal and Information
  Processing}, pages 945--948, 2013.

\bibitem{villani2009optimal}
C{\'e}dric Villani et~al.
\newblock {\em Optimal transport: old and new}, volume 338.
\newblock Springer, 2009.

\bibitem{vincent2011connection}
Pascal Vincent.
\newblock A connection between score matching and denoising autoencoders.
\newblock {\em Neural Computation}, 23(7):1661--1674, 2011.

\bibitem{wang2016learning}
Zhangyang Wang, Qing Ling, and Thomas Huang.
\newblock Learning deep l0 encoders.
\newblock {\em Proceedings of the AAAI Conference on Artificial Intelligence},
  30(1), 2016.

\bibitem{whang2021composing}
Jay Whang, Erik Lindgren, and Alex Dimakis.
\newblock Composing normalizing flows for inverse problems.
\newblock In Marina Meila and Tong Zhang, editors, {\em Proceedings of the 38th
  International Conference on Machine Learning}, volume 139 of {\em Proceedings
  of Machine Learning Research}, pages 11158--11169. PMLR, 2021.

\bibitem{winkler2019learning}
Christina Winkler, Daniel Worrall, Emiel Hoogeboom, and Max Welling.
\newblock Learning likelihoods with conditional normalizing flows.
\newblock {\em arXiv preprint arXiv:1912.00042}, 2019.

\bibitem{wu2019online}
Zihui Wu, Yu~Sun, Jiaming Liu, and Ulugbek Kamilov.
\newblock Online regularization by denoising with applications to phase
  retrieval.
\newblock In {\em Proceedings of the IEEE/CVF International Conference on
  Computer Vision (ICCV) Workshops}, 2019.

\bibitem{xiang2021fista}
Jinxi Xiang, Yonggui Dong, and Yunjie Yang.
\newblock Fista-net: Learning a fast iterative shrinkage thresholding network
  for inverse problems in imaging.
\newblock {\em IEEE Transactions on Medical Imaging}, 40(5):1329--1339, 2021.

\bibitem{yang2018dagan}
Guang Yang, Simiao Yu, Hao Dong, Greg Slabaugh, Pier~Luigi Dragotti, Xujiong
  Ye, Fangde Liu, Simon Arridge, Jennifer Keegan, Yike Guo, and David Firmin.
\newblock Dagan: Deep de-aliasing generative adversarial networks for fast
  compressed sensing {MRI} reconstruction.
\newblock {\em IEEE Transactions on Medical Imaging}, 37(6):1310--1321, 2018.

\bibitem{yang2018low}
Qingsong Yang, Pingkun Yan, Yanbo Zhang, Hengyong Yu, Yongyi Shi, Xuanqin Mou,
  Mannudeep~K. Kalra, Yi~Zhang, Ling Sun, and Ge~Wang.
\newblock Low-dose ct image denoising using a generative adversarial network
  with wasserstein distance and perceptual loss.
\newblock {\em IEEE Transactions on Medical Imaging}, 37(6):1348--1357, 2018.

\bibitem{yang2016deep}
Yan Yang, Jian Sun, Huibin Li, and Zongben Xu.
\newblock Deep admm-net for compressive sensing {MRI}.
\newblock In D.~Lee, M.~Sugiyama, U.~Luxburg, I.~Guyon, and R.~Garnett,
  editors, {\em Advances in Neural Information Processing Systems}, volume~29.
  Curran Associates, Inc., 2016.

\bibitem{yang2020admm}
Yan Yang, Jian Sun, Huibin Li, and Zongben Xu.
\newblock Admm-csnet: A deep learning approach for image compressive sensing.
\newblock {\em IEEE Transactions on Pattern Analysis and Machine Intelligence},
  42(3):521--538, 2020.

\bibitem{ye2018deep}
Dong~Hye Ye, Somesh Srivastava, Jean-Baptiste Thibault, Ken Sauer, and Charles
  Bouman.
\newblock Deep residual learning for model-based iterative ct reconstruction
  using plug-and-play framework.
\newblock In {\em 2018 IEEE International Conference on Acoustics, Speech and
  Signal Processing (ICASSP)}, pages 6668--6672, 2018.

\bibitem{yeh2017semantic}
Raymond~A. Yeh, Chen Chen, Teck Yian~Lim, Alexander~G. Schwing, Mark
  Hasegawa-Johnson, and Minh~N. Do.
\newblock Semantic image inpainting with deep generative models.
\newblock In {\em Proceedings of the IEEE Conference on Computer Vision and
  Pattern Recognition (CVPR)}, 2017.

\bibitem{yoo2021time}
Jaejun Yoo, Kyong~Hwan Jin, Harshit Gupta, Jerome Yerly, Matthias Stuber, and
  Michael Unser.
\newblock Time-dependent deep image prior for dynamic {MRI}.
\newblock {\em IEEE Transactions on Medical Imaging}, 40(12):3337--3348, 2021.

\bibitem{you2020ct}
Chenyu You, Guang Li, Yi~Zhang, Xiaoliu Zhang, Hongming Shan, Mengzhou Li,
  Shenghong Ju, Zhen Zhao, Zhuiyang Zhang, Wenxiang Cong, Michael~W. Vannier,
  Punam~K. Saha, Eric~A. Hoffman, and Ge~Wang.
\newblock Ct super-resolution gan constrained by the identical, residual, and
  cycle learning ensemble (gan-circle).
\newblock {\em IEEE Transactions on Medical Imaging}, 39(1):188--203, 2020.

\bibitem{zach2023stable}
Martin Zach, Florian Knoll, and Thomas Pock.
\newblock Stable deep {MRI} reconstruction using generative priors.
\newblock {\em IEEE Transactions on Medical Imaging}, pages 1--1, 2023.

\bibitem{zach2022computed}
Martin Zach, Erich Kobler, and Thomas Pock.
\newblock Computed tomography reconstruction using generative energy-based
  priors.
\newblock {\em arXiv preprint arXiv:2203.12658}, 2022.

\bibitem{zach2023explicit}
Martin Zach, Thomas Pock, Erich Kobler, and Antonin Chambolle.
\newblock Explicit diffusion of gaussian mixture model based image priors.
\newblock In Luca Calatroni, Marco Donatelli, Serena Morigi, Marco Prato, and
  Matteo Santacesaria, editors, {\em Scale Space and Variational Methods in
  Computer Vision}, pages 3--15, Cham, 2023. Springer International Publishing.

\bibitem{zhang2020review}
Hai-Miao Zhang and Bin Dong.
\newblock A review on deep learning in medical image reconstruction.
\newblock {\em Journal of the Operations Research Society of China},
  8:311--340, 2020.

\bibitem{zhang2018ista}
Jian Zhang and Bernard Ghanem.
\newblock Ista-net: Interpretable optimization-inspired deep network for image
  compressive sensing.
\newblock In {\em Proceedings of the IEEE Conference on Computer Vision and
  Pattern Recognition (CVPR)}, 2018.

\bibitem{zhang2022plug}
Kai Zhang, Yawei Li, Wangmeng Zuo, Lei Zhang, Luc Van~Gool, and Radu Timofte.
\newblock Plug-and-play image restoration with deep denoiser prior.
\newblock {\em IEEE Transactions on Pattern Analysis and Machine Intelligence},
  44(10):6360--6376, 2022.

\bibitem{zhang2017beyond}
Kai Zhang, Wangmeng Zuo, Yunjin Chen, Deyu Meng, and Lei Zhang.
\newblock Beyond a gaussian denoiser: Residual learning of deep cnn for image
  denoising.
\newblock {\em IEEE Transactions on Image Processing}, 26(7):3142--3155, 2017.

\bibitem{zhang2017learning}
Kai Zhang, Wangmeng Zuo, Shuhang Gu, and Lei Zhang.
\newblock Learning deep cnn denoiser prior for image restoration.
\newblock In {\em Proceedings of the IEEE Conference on Computer Vision and
  Pattern Recognition (CVPR)}, 2017.

\bibitem{zhang2021residual}
Yulun Zhang, Yapeng Tian, Yu~Kong, Bineng Zhong, and Yun Fu.
\newblock Residual dense network for image restoration.
\newblock {\em IEEE Transactions on Pattern Analysis and Machine Intelligence},
  43(7):2480--2495, 2021.

\bibitem{thu2018image}
Bo~Zhu, Jeremiah~Z. Liu, Stephen~F. Cauley, Bruce~R. Rosen, and Matthew~S.
  Rosen.
\newblock Image reconstruction by domain-transform manifold learning.
\newblock {\em Nature}, 555(7697):487--492, 2018.

\bibitem{zoran2011from}
Daniel Zoran and Yair Weiss.
\newblock From learning models of natural image patches to whole image
  restoration.
\newblock In {\em 2011 International Conference on Computer Vision}, pages
  479--486, 2011.

\end{thebibliography}
